\documentclass[11pt, a4paper]{article}
\usepackage[cp1251]{inputenc}
\usepackage{amsmath} \usepackage{euscript} \usepackage{marvosym}
\usepackage[dvipsnames]{xcolor}
\usepackage{mymatrix2}%\let\mymatrixII\mymatrix
\usepackage{graphicx}
\usepackage{longtable}
\usepackage{mathrsfs}
\usepackage{amssymb, amscd}
\usepackage{comment}
\begin{document}

\newcounter{bnomer} \newcounter{snomer}
\newcounter{bsnomer}
\setcounter{bnomer}{0}
\renewcommand{\thesnomer}{\thebnomer.\arabic{snomer}}
\renewcommand{\thebsnomer}{\thebnomer.\arabic{bsnomer}}
\renewcommand{\refname}{\begin{center}\large{\textbf{References}}\end{center}}

\setcounter{MaxMatrixCols}{14}

\newcommand\restr[2]{{% we make the whole thing an ordinary symbol
  \left.\kern-\nulldelimiterspace % automatically resize the bar with \right
  #1 % the function
  %\vphantom{\big|} % pretend it's a little taller at normal size
  \right|_{#2} % this is the delimiter
}}

\newcommand{\sect}[1]{%
\setcounter{snomer}{0}\setcounter{bsnomer}{0}
\refstepcounter{bnomer}
\par\bigskip\begin{center}\large{\textbf{\arabic{bnomer}. {#1}}}\end{center}}
\newcommand{\sst}[1]{%
\refstepcounter{bsnomer}
\par\bigskip\textbf{\arabic{bnomer}.\arabic{bsnomer}. {#1}}\par}
\newcommand{\defi}[1]{%
\refstepcounter{snomer}
\par\medskip\textbf{Definition \arabic{bnomer}.\arabic{snomer}. }{#1}\par\medskip}
\newcommand{\theo}[2]{%
\refstepcounter{snomer}
\par\textbf{Theorem \arabic{bnomer}.\arabic{snomer}. }{#2} {\emph{#1}}\hspace{\fill}$\square$\par}
\newcommand{\mtheop}[2]{%
\refstepcounter{snomer}
\par\textbf{Theorem \arabic{bnomer}.\arabic{snomer}. }{\emph{#1}}
\par\textsc{Proof}. {#2}\hspace{\fill}$\square$\par}
\newcommand{\mcorop}[2]{%
\refstepcounter{snomer}
\par\textbf{Corollary \arabic{bnomer}.\arabic{snomer}. }{\emph{#1}}
\par\textsc{Proof}. {#2}\hspace{\fill}$\square$\par}
\newcommand{\mtheo}[1]{%
\refstepcounter{snomer}
\par\medskip\textbf{Theorem \arabic{bnomer}.\arabic{snomer}. }{\emph{#1}}\par\medskip}
\newcommand{\theobn}[1]{%
\par\medskip\textbf{Theorem. }{\emph{#1}}\par\medskip}
\newcommand{\theoc}[2]{%
\refstepcounter{snomer}
\par\medskip\textbf{Theorem \arabic{bnomer}.\arabic{snomer}. }{#1} {\emph{#2}}\par\medskip}
\newcommand{\mlemm}[1]{%
\refstepcounter{snomer}
\par\medskip\textbf{Lemma \arabic{bnomer}.\arabic{snomer}. }{\emph{#1}}\par\medskip}
\newcommand{\mprop}[1]{%
\refstepcounter{snomer}
\par\medskip\textbf{Proposition \arabic{bnomer}.\arabic{snomer}. }{\emph{#1}}\par\medskip}
\newcommand{\theobp}[2]{%
\refstepcounter{snomer}
\par\textbf{Theorem \arabic{bnomer}.\arabic{snomer}. }{#2} {\emph{#1}}\par}
\newcommand{\theop}[2]{%
\refstepcounter{snomer}
\par\textbf{Theorem \arabic{bnomer}.\arabic{snomer}. }{\emph{#1}}
\par\textsc{Proof}. {#2}\hspace{\fill}$\square$\par}
\newcommand{\theosp}[2]{%
\refstepcounter{snomer}
\par\textbf{Theorem \arabic{bnomer}.\arabic{snomer}. }{\emph{#1}}
\par\textsc{Sketch of the proof}. {#2}\hspace{\fill}$\square$\par}
\newcommand{\exam}[1]{%
\refstepcounter{snomer}
\par\medskip\textbf{Example \arabic{bnomer}.\arabic{snomer}. }{#1}\par\medskip}
\newcommand{\deno}[1]{%
\refstepcounter{snomer}
\par\textbf{Notation \arabic{bnomer}.\arabic{snomer}. }{#1}\par}
\newcommand{\lemm}[1]{%
\refstepcounter{snomer}
\par\textbf{Lemma \arabic{bnomer}.\arabic{snomer}. }{\emph{#1}}\hspace{\fill}$\square$\par}
\newcommand{\lemmp}[2]{%
\refstepcounter{snomer}
\par\medskip\textbf{Lemma \arabic{bnomer}.\arabic{snomer}. }{\emph{#1}}
\par\textsc{Proof}. {#2}\hspace{\fill}$\square$\par\medskip}
\newcommand{\coro}[1]{%
\refstepcounter{snomer}
\par\textbf{Corollary \arabic{bnomer}.\arabic{snomer}. }{\emph{#1}}\hspace{\fill}$\square$\par}
\newcommand{\mcoro}[1]{%
\refstepcounter{snomer}
\par\textbf{Corollary \arabic{bnomer}.\arabic{snomer}. }{\emph{#1}}\par\medskip}
\newcommand{\corop}[2]{%
\refstepcounter{snomer}
\par\textbf{Corollary \arabic{bnomer}.\arabic{snomer}. }{\emph{#1}}
\par\textsc{Proof}. {#2}\hspace{\fill}$\square$\par}
\newcommand{\nota}[1]{%
\refstepcounter{snomer}
\par\medskip\textbf{Remark \arabic{bnomer}.\arabic{snomer}. }{#1}\par\medskip}
\newcommand{\propp}[2]{%
\refstepcounter{snomer}
\par\medskip\textbf{Proposition \arabic{bnomer}.\arabic{snomer}. }{\emph{#1}}
\par\textsc{Proof}. {#2}\hspace{\fill}$\square$\par\medskip}
\newcommand{\hypo}[1]{%
\refstepcounter{snomer}
\par\medskip\textbf{Conjecture \arabic{bnomer}.\arabic{snomer}. }{\emph{#1}}\par\medskip}
\newcommand{\prop}[1]{%
\refstepcounter{snomer}
\par\textbf{Proposition \arabic{bnomer}.\arabic{snomer}. }{\emph{#1}}\hspace{\fill}$\square$\par}

\newcommand{\proof}[2]{%
\par\medskip\textsc{Proof{#1}}. \hspace{-0.2cm}{#2}\hspace{\fill}$\square$\par\medskip}

\makeatletter
\def\iddots{\mathinner{\mkern1mu\raise\p@
\vbox{\kern7\p@\hbox{.}}\mkern2mu
\raise4\p@\hbox{.}\mkern2mu\raise7\p@\hbox{.}\mkern1mu}}
\makeatother

\newcommand{\okr}[2]{%
\refstepcounter{snomer}
\par\medskip\textbf{{#1} \arabic{bnomer}.\arabic{snomer}. }{\emph{#2}}\par\medskip}

\newcommand{\Ind}[3]{%
\mathrm{Ind}_{#1}^{#2}{#3}}
\newcommand{\Res}[3]{%
\mathrm{Res}_{#1}^{#2}{#3}}
\newcommand{\epsi}{\varepsilon}
\newcommand{\tri}{\triangleleft}
\newcommand{\Supp}[1]{%
\mathrm{Supp}(#1)}
\newcommand{\SSu}[1]{%
\mathrm{SingSupp}(#1)}

\newcommand{\gee}{\geqslant}
\newcommand{\reg}{\mathrm{reg}}
\newcommand{\Dyn}{\mathrm{Dyn}}
\newcommand{\Ann}{\mathrm{Ann}\,}
\newcommand{\Cent}[1]{\mathbin\mathrm{Cent}({#1})}
\newcommand{\PCent}[1]{\mathbin\mathrm{PCent}({#1})}
\newcommand{\Irr}[1]{\mathbin\mathrm{Irr}({#1})}
\newcommand{\Exp}[1]{\mathbin\mathrm{Exp}({#1})}
\newcommand{\empr}[2]{[-{#1},{#1}]\times[-{#2},{#2}]}
\newcommand{\sreg}{\mathrm{sreg}}
\newcommand{\ilm}{\varinjlim}
\newcommand{\wdth}{\mathrm{wd}}
\newcommand{\plm}{\varprojlim}
\newcommand{\codim}{\mathrm{codim}\,}
\newcommand{\GKdim}{\mathrm{GKdim}\,}
\newcommand{\chara}{\mathrm{char}\,}
\newcommand{\rk}{\mathrm{rk}\,}
\newcommand{\chr}{\mathrm{ch}\,}
\newcommand{\Ker}{\mathrm{Ker}\,}
\newcommand{\id}{\mathrm{id}}
\newcommand{\Ad}{\mathrm{Ad}}
\newcommand{\Gh}{\mathrm{Gh}}
\newcommand{\col}{\mathrm{col}}
\newcommand{\row}{\mathrm{row}}
\newcommand{\high}{\mathrm{high}}
\newcommand{\low}{\mathrm{low}}
\newcommand{\pho}{\hphantom{\quad}\vphantom{\mid}}
\newcommand{\fho}[1]{\vphantom{\mid}\setbox0\hbox{00}\hbox to \wd0{\hss\ensuremath{#1}\hss}}
\newcommand{\wt}{\widetilde}
\newcommand{\wh}{\widehat}
\newcommand{\ad}[1]{\mathrm{ad}_{#1}}
\newcommand{\tr}{\mathrm{tr}\,}
\newcommand{\GL}{\mathrm{GL}}
\newcommand{\SL}{\mathrm{SL}}
\newcommand{\SO}{\mathrm{SO}}
\newcommand{\Or}{\mathrm{O}}
\newcommand{\Sp}{\mathrm{Sp}}
\newcommand{\Sa}{\mathrm{S}}
\newcommand{\Ua}{\mathrm{U}}
\newcommand{\Andre}{\mathrm{Andre}}
\newcommand{\Aord}{\mathrm{Aord}}
\newcommand{\Mat}{\mathrm{Mat}}
\newcommand{\Stab}{\mathrm{Stab}}
\newcommand{\htt}{\mathfrak{h}}
\newcommand{\spt}{\mathfrak{sp}}
\newcommand{\slt}{\mathfrak{sl}}
\newcommand{\sot}{\mathfrak{so}}

\newcommand{\vfi}{\varphi}
\newcommand{\aad}{\mathrm{ad}}
\newcommand{\vpi}{\varpi}
\newcommand{\teta}{\vartheta}
\newcommand{\Bfi}{\Phi}
\newcommand{\Fp}{\mathbb{F}}
\newcommand{\Rp}{\mathbb{R}}
\newcommand{\Zp}{\mathbb{Z}}
\newcommand{\Cp}{\mathbb{C}}
\newcommand{\Ap}{\mathbb{A}}
\newcommand{\Pp}{\mathbb{P}}
\newcommand{\Kp}{\mathbb{K}}
\newcommand{\Np}{\mathbb{N}}
\newcommand{\ut}{\mathfrak{u}}
\newcommand{\at}{\mathfrak{a}}
\newcommand{\glt}{\mathfrak{gl}}
\newcommand{\hei}{\mathfrak{hei}}
\newcommand{\nt}{\mathfrak{n}}
\newcommand{\kt}{\mathfrak{k}}
\newcommand{\mt}{\mathfrak{m}}
\newcommand{\rt}{\mathfrak{r}}
\newcommand{\rad}{\mathfrak{rad}}
\newcommand{\bt}{\mathfrak{b}}
\newcommand{\unt}{\underline{\mathfrak{n}}}
\newcommand{\gt}{\mathfrak{g}}
\newcommand{\vt}{\mathfrak{v}}
\newcommand{\pt}{\mathfrak{p}}
\newcommand{\Xt}{\mathfrak{X}}
\newcommand{\Po}{\mathcal{P}}
\newcommand{\PV}{\mathcal{PV}}
\newcommand{\Uo}{\EuScript{U}}
\newcommand{\Fo}{\EuScript{F}}
\newcommand{\Do}{\EuScript{D}}
\newcommand{\Eo}{\EuScript{E}}
\newcommand{\Jo}{\EuScript{J}}
\newcommand{\Iu}{\mathcal{I}}
\newcommand{\Mo}{\mathcal{M}}
\newcommand{\Nu}{\mathcal{N}}
\newcommand{\Ro}{\mathcal{R}}
\newcommand{\Co}{\mathcal{C}}
\newcommand{\Ko}{\mathcal{K}}
\newcommand{\So}{\mathcal{S}}
\newcommand{\Lo}{\mathcal{L}}
\newcommand{\Ou}{\mathcal{O}}
\newcommand{\Uu}{\mathcal{U}}
\newcommand{\Tu}{\mathcal{T}}
\newcommand{\Au}{\mathcal{A}}
\newcommand{\Vu}{\mathcal{V}}
\newcommand{\Du}{\mathcal{D}}
\newcommand{\Bu}{\mathcal{B}}
\newcommand{\Sy}{\mathcal{Z}}
\newcommand{\Sb}{\mathcal{F}}
\newcommand{\Gr}{\mathcal{G}}
\newcommand{\Xu}{\mathcal{X}}
\newcommand{\Op}{\mathbb{O}}
\newcommand{\chv}{\mathrm{chv}}
\newcommand{\rtc}[1]{C_{#1}^{\mathrm{red}}}

\author{Mikhail Ignatev\and Mikhail Venchakov}
\date{}
\title{Orbits and characters associated with rook placements\\ for Sylow $p$-subgroups of finite orthogonal groups}\maketitle
\begin{abstract} Let $U$ be a Sylow $p$-subgroup in a classical group over a finite field of characteristic $p$. The coadjoint orbits of the group $U$ play the key role in the description of irreducible complex characters of $U$. Almost all important classes of orbits and characters studied to the moment can be uniformly described as the orbits and characters associated with so-called orthogonal rook placements. In the paper, we study such orbits for the orthogonal group. We construct a polarization for the canonical form on such an orbit and present a semi-direct decomposition for the corresponding irreducible characters in the spirit of the Mackey little group method. As a corollary, we compute the dimension of an orbit associated with an orthogonal rook placement.

\end{abstract}

\sect{Introduction}

\let\thefootnote\relax\footnote{The research was supported by RSF (project No. 25--21--00219), \texttt{https://rscf.ru/en/project/25-21-00219/}.}

Let $U$ be a unipotent algebraic group over a finite field $\Fp_q$ of sufficiently large characteristic $p$. The main tool in representation theory of $U$ is the orbit method created in 1962 by A.A. Kirillov, see~\cite{Kirillov62}, \cite{Kirillov04}, \cite{Kazhdan77}. The key idea of the orbit method says that the irreducible representations of $U$ are in one-to-one correspondence with the coadjoint orbits of this group. Namely, the group $U$ acts on its Lie algebra $\ut$ via the adjoint action; the dual action of $U$ on the dual space $\ut^*$ is called coadjoint. It turns out that there is a natural bijection between the set $\Irr{U}$ of all irreducible complex characters of $U$ and the set $\ut^*/U$ of coadjoint orbits, see Theorem~\ref{theo:orbit_method} below for the precise statement.

Let $U$ be a maximal unipotent subgroup (or, equivalently, a Sylow $p$-subgroup) in a simple classical group over $\Fp_q$. A complete description of all coadjoint orbits for the group $U$ is a wild problem, so a natural question is how to describe certain important classes of orbits and the corresponding irreducible characters. For the case $A_{n-1}$, orbits of maximal possible dimension were classified in the first Kirillov's work on the orbit method \cite{Kirillov62}. All of them are associated with the Kostant cascade, which is in fact an orthogonal rook placement, see the next section for the precise definitions. For other classical root systems, orbits associated with the Kostant cascades also have maximal possible dimension, see \cite{Kostant12} and \cite{Kostant13}. For $A_{n-1}$, orbits of submaximal dimension were classified by A.N. Panov \cite{IgnatevPanov09}. Such orbits also correspond to orthogonal rook placements (modulo adding simple root covectors to the canonical forms on them). A classification of orbits of maximal dimension in type $C_n$ follows from C. Andre and A. Neto's paper \cite{AndreNeto06}. A description of orbits of maximal possible dimension for $B_n$ and $D_n$, as well as of orbits of submaximal dimension for all types $B_n$, $C_n$ and $D_n$, will be presented in \cite{IgnatevPetukhovVenchakov24}. All these orbits are also associated with orthogonal rook placements (modulo adding simple root covectors to the canonical forms on them).

The irreducible characters corresponding to the orbits of maximal dimension $M$ in type $A_{n-1}$ were computed by C. Andre in 2011 \cite{Andre01}. The characters of submaximal dimension $M-2$ for $A_{n-1}$ were calculated by the first author in \cite{Ignatev09}. A formula for the characters of maximal dimension in type $C_n$ can be obtained from the results of Andre and Neto \cite{AndreNeto06}. Characters of maximal dimension in the orthogonal case, as well as the characters of submaximal dimension for classical cases, will be computed in \cite{IgnatevPetukhovVenchakov24}. Andre's method is based on very nice stratification of $\ut^*$ established by him, whose strata are numerated by arbitrary (possibly, non-orthogonal) rook placements and maps from them to the set $\Fp_q^{\times}$ of non-zero numbers.

In this paper, we present another approach to calculating characters based on the Mackey method of little groups. This method reduces the calculation of irreducible characters of a semi-direct product group with abelian normal subgroup to the calculation of irreducible characters of certain subgroups called little groups, see Section~\ref{orbit_semidirect} for the detail. It turned out that this method allows to compute characters of maximal and submaximal dimension for $A_{n-1}$. Furthermore, in \cite{IgnatevVenchakov24} we applied this method to the calculation of characters of the next possible dimension $M-4$, while in \cite{IgnatevPetukhovVenchakov24} this method allows to obtain explicit formulas for characters of maximal and submaximal dimension in the orthogonal and symplectic cases. Here we apply the Mackey method to the characters associated with arbitrary orthogonal rook placements in the orthogonal case.

The structure of the paper is as follows. In Section~\ref{main_definitions}, we briefly recall main definitions and facts about finite orthogonal groups and its maximal unipotent subgroups. We also define rook placements in root systems and associated coadjoint orbits. In Section~\ref{sect:sea_battle_pol}, we recall the orbit method. In particular, to construct the representation corresponding to the orbit of a linear form $f$, it is needed to present a polarization for $f$. Our first main result established such a polarization for the canonical form on an orbit associated with an orthogonal rook placement, see Theorem~\ref{theo_pol}. As a corollary, we obtain a formula for the dimension of the orbit of $f$, see Corollary~\ref{coro_dim_pol}.

Section~\ref{orbit_semidirect} is central. First, we briefly recall the Mackey method of little groups for semi-direct products with abelian normal subgroup, see Theorem~\ref{theo_semi_direct}. Our second main result, Theorem~\ref{theo_char_decompose}, establishes a decomposition in the spirit of the Mackey method for the irreducible character corresponding to an orbit associated with an orthogonal rook placement. As an application, in Section~\ref{formula_dimension} we present another explicit formula for the dimension of such an orbit. The answer is given in terms of the involution equal to the product of all reflections corresponding to the roots from this rook placement, see Theorem~\ref{theo_dim_orbit}.

\sect{Main definitions} \label{main_definitions}
To begin with we present some basic facts about root systems of orthogonal algebras. We will denote root systems of type $B_n$ or $D_n$ by $\Phi$. As usual, we will identify them with the following subsets of $\Rp^n$:
$$
\begin{array}{ll}
&B_n = \{\pm\epsi_i\pm\epsi_j, 1\leq i < j\leq
n\}\cup\{\pm\epsi_i, 1\leq i\leq n\},\\
&D_n = \{\pm\epsi_i\pm\epsi_j, 1\leq i < j\leq n\},
\end{array}
$$
where $\{\epsi_i\}_{i=1}^n$ is the standard basis in $\Rp^n$. Pick the set of simple roots $\Delta=\Delta(\Phi)$ as in
\cite{Bourbaki03}:
$$
\begin{array}{ll}
&\Delta(B_n)=\{\epsi_i-\epsi_{i+1},1\leq i\leq n-1\}\cup\{\epsi_n\},\\
&\Delta(D_n)=\{\epsi_i-\epsi_{i+1},1\leq i\leq n-1\}\cup\{\epsi_{n-1}+\epsi_n\}.\\
\end{array}
$$
The set of positive roots $\Phi^+\supset\Delta$ looks as follows:
$$
\begin{array}{ll}
&B_n^+ = \{\epsi_i\pm\epsi_j, 1\leq i < j\leq n\}\cup\{\epsi_i,
1\leq i\leq n\},\\
&D_n^+ = \{\epsi_i\pm\epsi_j, 1\leq i < j\leq n\}.\\
\end{array}
$$
Let
$$
m=\left\{\begin{array}{ll}2n+1,&\mbox{if }\Phi=B_n,\\
2n,&\mbox{if }\Phi=D_n.
\end{array}\right.
$$
Denote by $\ut=\ut(\Phi)$ the subalgebra of $\mathfrak{gl}_m(\Fp_q)$
spanned by the vectors $e_{\alpha}$, $\alpha\in\Phi$, where
$$
\begin{array}{ll}
&e_{\epsi_i}=e_{0, i}-e_{-i,0},\quad 1\leq i\leq n,\\
&e_{\epsi_i-\epsi_j}=e_{j,i}-e_{-i,-j},\quad 1\leq i<j\leq n,\\
&e_{\epsi_i+\epsi_j}=e_{-j,i}-e_{-i,j},\quad 1\leq i<j\leq n.\\
\end{array}
$$
Here we numerate the rows and the columns of $m\times m$ matrix by the indices $$1,2,\ldots,n,0,-n,\ldots,-2,-1$$ (there is no index $0$ in the case $D_n$), and we denote by $e_{a, b}$ the usual elementary matrix. Hence, $\ut$ is the algebra of lower triangular matrices with zeroes on the diagonal that are skew-symmetric with respect to the antidiagonal. Of course, it is a maximal nilpotent subalgebra in the corresponding classical algebra $\gt= \gt(\Phi)$. In particular, $\dim\ut = |\Phi^+|$. \exam{Here we schematically drew the algebras $B_3$ and $D_4$. The convention is as follows: given $1\leq j<i$, the square $(i,j)$ (respectively, $(-i,j)$ and $(0,j)$) corresponds to the root $\epsi_j-\epsi_i$ (respectively, $\epsi_j+\epsi_i$ and $\epsi_j$).
\begin{center}\small
$\mymatrix{
\pho& \pho&\pho& \pho& \pho& \pho & \pho\\
\Top{2pt}\Rt{2pt} \pho & \pho& \pho& \pho& \pho& \pho & \pho\\
\pho & \Top{2pt}\Rt{2pt} \pho& \pho& \pho& \pho& \pho & \pho\\
\pho & \pho& \Top{2pt}\Rt{2pt} \pho& \pho& \pho& \pho & \pho\\
\pho & \pho& \Top{2pt}\Lft{2pt} 0& \Top{2pt}\Rt{2pt} \pho& \pho& \pho & \pho\\
\pho& \Top{2pt}\Lft{2pt} 0& \pho& \pho& \Top{2pt}\Rt{2pt} \pho& \pho & \pho\\
\Top{2pt} 0& \pho& \pho& \pho& \pho& \Top{2pt}\Rt{2pt}\pho & \pho\\
}\quad \mymatrix{
\pho& \pho&\pho& \pho& \pho& \pho & \pho & \pho\\
\Top{2pt}\Rt{2pt} \pho & \pho& \pho& \pho& \pho& \pho & \pho & \pho\\
\pho & \Top{2pt}\Rt{2pt} \pho& \pho& \pho& \pho& \pho & \pho & \pho\\
\pho & \pho& \Top{2pt}\Rt{2pt} \pho& \pho& \pho& \pho & \pho & \pho\\
\pho &\pho & \pho& \Top{2pt}\Lft{2pt}\Rt{2pt} 0& \pho& \pho& \pho & \pho\\
\pho &\pho& \Top{2pt}\Lft{2pt} 0& \pho& \Top{2pt}\Rt{2pt}\pho& \pho& \pho & \pho\\
\pho &\Top{2pt}\Lft{2pt} 0& \pho& \pho& \pho& \Top{2pt}\Rt{2pt} \pho& \pho & \pho\\
\Top{2pt} 0& \pho& \pho& \pho& \pho& \pho & \Top{2pt}\Rt{2pt}\pho &\pho\\
}$
\end{center}
}

Furthermore, we define the functions 
$$
\begin{array}{ll}
&\col\colon\Phi^+\to\{1,\ldots,n\}\colon\col(\epsi_i\pm\epsi_j)=\col(\epsi_i)=i,\\
&\row\colon\Phi^+\to\{-n,\ldots,n\}\colon\row(\epsi_i\pm\epsi_j)=\mp j,\row(\epsi_i)=0.\\
\end{array}
$$
For arbitrary $-n+1\leq i\leq n-1$ and $1\leq j\leq n$, the sets
$$
\begin{array}{ll}
&R_i = R_i(\Phi) = \{\alpha\in\Phi^+\mid \row(\alpha)=i\},\\
&C_j = C_j(\Phi) = \{\alpha\in\Phi^+\mid \col(\alpha)=j\}\\
\end{array}
$$
are called the $i$th\emph{ row} and the $j$th\emph{ column}
$\Phi^+$ respectively. We introduce the mirror order on the 
set of indices $$1\prec2\prec\ldots\prec
n\prec0\prec-n\prec\ldots\prec-2\prec-1,$$ and the following total orders on $\Phi^+$:
\begin{equation}
\begin{array}{ll}
&\alpha\prec\beta~\stackrel{\mathrm{def}}{\iff}~
\col(\beta)\prec\col(\alpha)~\mbox{or}~
\col(\beta)=\col(\alpha),~\row(\beta)\prec\row(\alpha),\\
&\alpha\prec'\beta~\stackrel{\mathrm{def}}{\iff}~
\col(\beta)\prec\col(\alpha)~\mbox{or}~
\col(\beta)=\col(\alpha),~\row(\beta)\succ\row(\alpha).\label{formula_complete_orders}
\end{array}
\end{equation}
For example, for $\Phi=B_6$ we have
$\epsi_2-\epsi_4\succ\epsi_2\succ\epsi_2+\epsi_5\succ\epsi_3-\epsi_6$,
$\epsi_2+\epsi_5\succ'\epsi_2\succ'\epsi_2-\epsi_4\succ'\epsi_3-\epsi_6$.

Now we are ready to define the main object of our interest. \defi{Let a subset $D=\{\beta_1,\ldots,\beta_t\}\subset\Phi^+$ consisting of pairwise orthogonal roots satisfy the following condition:
\begin{equation}
|D\cap R_i|\leq 1\mbox{ and }|D\cap C_j|\leq 1\mbox{ for all }i,
j.\label{formula_basic_supp}
\end{equation}
We call such a subset $D$ an \emph{orthogonal rook placement}. (Note than in \cite{AndreNeto06}, such a subset $D$ is called a basic subset of $\Phi$.)}

\nota{Let $W=W(\Phi)$ be the Weyl group of the root system $\Phi$. To each orthogonal rook placement $D=\{\beta_1,\ldots,\beta_t\}$ one can attach the involution (i.e., the element of order 2)
$\sigma=\sigma_D\in W$:
\begin{equation}
\sigma=r_{\beta_1}\ldots r_{\beta_t}.
\label{formula_ortog_decompose}
\end{equation}
Here, for an arbitrary $\beta\in\Phi$, we denote by $r_{
\beta}$ the reflection of $\Rp^n$ with respect to the hyperplane orthogonal to  $\beta$. We will call the set $D$ the \emph{support} of the  involution $\sigma$ and write $D=\Supp{\sigma}$.

%fix
%Note that, in the case of $A_n$ or $C_n$, each involution can be obtained in such a way. On the other hand, for $B_n$ or $D_n$, this is not true. For example, if $\Phi=B_2$ and $\sigma=r_{\epsi_1}\cdot r_{\epsi_2}$, then there is no orthogonal rook placement $D$ such that $D=\Supp{\sigma}$.

Note that, in the case of $B_n$ or $D_n$, not each involution can be obtained in such a way. For example, if $\Phi=B_2$ and $\sigma=r_{\epsi_1}\cdot r_{\epsi_2}$, then there is no orthogonal rook placement $D$ such that $D=\Supp{\sigma}$.
}

Denote the maximal unipotent subgroup $\exp(\ut)$ in the corresponding classical finite group $G$ by~$U$. In the sequel, we will assume everywhere that $p=\chara{\Fp_q}\gee n$. Under this assumption, the map
$$
\exp\colon\ut\to U\colon x\mapsto \sum_{i=0}^{m-1}\frac{x^i}{i!}
$$
is well-defined and is in fact a bijection (and also an isomorphism of algebraic varieties over $\overline{\Fp_q}$);
we denote the inverse map by $\ln$. Furthermore, the 
Backer--Campbell--Hausdorff formula claims that, for a Lie subalgebra $\at\subset\ut$ and arbitrary $u, v\in\at$, one has
\begin{equation}
\exp(u)\exp(v)=\exp(u + v + \tau(u,
v)),\label{formula_Campbell_Hausdorf}
\end{equation} where
$\tau(u, v)\in[\at, \at]$ (here $[\at, \at]=\langle[x, y], x,
y\in\at\rangle_{\Fp_q}$).

The group $U$ acts on its Lie algebra $\ut$ by the adjoint action; the dual action of $U$ on the $\Fp_q$-dual space $\ut^*$ is called \emph{coadjoint}. Using the non-degenerate form $$\langle A, B\rangle=\frac{1}{2}\mathrm{tr}(AB)$$ on $\mathfrak{gl}_m(\Fp_q)$, one can identify the dual space $\ut^*$ with the space $\ut^t$ (in this case, $e_{\alpha}^*=e_{\alpha}^t$ for any
$\alpha\in\Phi^+$). Under this identification, the coadjoint action has the following form:
$$
g.x = \mathrm{pr}(gxg^{-1}),~g\in U,~x\in\ut^*
$$
(here we denote the projection
$\mathfrak{gl}_m(\Fp_q)\to\ut^*$ along $\ut$ by $\mathrm{pr}$). The orbits of the coadjoint action play the key role in the classification of irreducible representations of group $U$ (see Section~\ref{sect:sea_battle_pol}).
\defi{Let $D=\{\beta_1,\ldots,\beta_t\}\subset\Phi^+$ be an orthogonal rook placement, and
$\xi=(\xi_{\beta})_{\beta\in D}$ be a set of non-zero constants from $\Fp_q$. Put
$$
f=f_{D, \xi} = \sum_{\beta\in D}\xi_{\beta}e_{\beta}^*\in\ut^*. $$ We say that the coadjoint orbit $\Omega=\Omega_{D,
\xi}\subset\ut^*$ of the linear form $f$ is \emph{associated} with $D$, and $f$ is called the
\emph{canonical form} on $\Omega$.}

\sect{Battleship and polarizations} \label{sect:sea_battle_pol}

First, we will recall the main idea of the orbit method created by A.A. Kirillov \cite{Kirillov62}, \cite{Kirillov04} and adapted by D. Kazhdan \cite{Kazhdan77} to finite groups (see also \cite{Kirillov95}), and construct a polarization for orbits associated with orthogonal rook placements.

Recall the main definitions. Let $\gt$ be an arbitrary Lie algebra, and $f\in\gt^*$ be an arbitrary linear form. A subalgebra
$\pt\subset\gt$ is called a \emph{polarization} for $f$ if $f([x,
y])=0$ for all $x$ , $y\in\pt$ (such subspaces are called $f$-\emph{isotropic}), and $\pt$ is a maximal $f$-isotropic subspace. It is easy to see that the dimension of the coadjoint orbit of $f$ is connected with the dimension of a polarization of $f$ by the simple formula $\dim\Omega_f=2\cdot\codim\pt$ (see, e.g., \cite[page
117]{Srinivasan}). According to the orbit method, the construction of a polarization is a necessary step to obtain the irreducible representation corresponding to a given linear form.

Specifically, pick and fix $\theta\colon\Fp_q\to\Cp^{\times}=\Cp\setminus\{0\}$, an arbitrary nontrivial character of the field $\Fp_q$ (i.e., a nontrivial homomorphism from the additive group of this field to the multiplicative group of complex numbers). Recall that we called the map inverse to $\exp\colon\ut\to U$ by $\ln\colon
U\to\ut$.
Let $\Omega\subset\ut^*$ be an arbitrary coadjoint orbit.
Consider the following function $\chi=\chi_{\Omega}\colon U\to\Cp$:
\begin{equation}
\chi(g)=q^{-\frac{1}{2}\dim\Omega}\cdot\sum_{f\in\Omega}\theta(f(\ln(g))),\quad
g\in U\label{formula_char_orbit}
\end{equation}
(over $\overline{\Fp}_q$, coadjoint orbits are affine varieties of even dimension in $\ut^*$; furthermore,
$q^{\dim\Omega}=|\Omega|$).
The main idea of the orbit method can be expressed by the following theorem. \theo{The map\label{theo:orbit_method}
$\Omega\mapsto\chi_{\Omega}$ establishes a one-to-one correspondence between the set of coadjoint orbits of the group $U$ and the set $\Irr{U}$ of finite-dimensional irreducible complex characters of the group $U$. Furthermore\textup, the complex dimension of the representation corresponding to $\Omega$ equals $q^{\frac{1}{2}\dim\Omega}=\sqrt{|\Omega|}.$}{\cite[Proposition
2]{Kazhdan77}} \label{theo_char_orbit} Moreover, the construction of the representation corresponding to a given orbit $\Omega$ is as follows. One should pick a linear form $f\in\Omega$ and a polarization $\pt$ of $\Omega$. Next, one should consider the following one-dimensional representation of $\exp(\pt)$:
$$\vfi_f(h)=\theta(f(\ln(h))),~h\in\exp(\pt).$$
The required representation now has the form $\Ind{\exp(\pt)}{U}(\vfi_f)$.

Here we present a construction of the polarization for the canonical form on an arbitrary orbit in $\ut^*$ associated with an orthogonal rook placement.
\defi{Let $\beta\in\Phi^+$. Roots $\alpha$, $\gamma\in\Phi^+$ are called
$\beta$-\emph{singular} if $\alpha+\gamma=\beta$. The set of all $\beta$-singular roots is denoted by $S(\beta)$ (see
\cite{Andre95ii} for $A_n$ and \cite{AndreNeto06} for
$B_n$, $D_n$).} It is easy to see that the singular roots look as follows:
\begin{equation}
\begin{array}{ll}
&S(\epsi_i-\epsi_j)=\bigcup_{l=i+1}^{j-1}\{\epsi_i-\epsi_l,\epsi_l-\epsi_j\},~1\leq i < j\leq n,\\
&S(\epsi_i) = \bigcup_{l=i+1}^{n}\{\epsi_i-\epsi_l,\epsi_l\},~1\leq i\leq n,\\
%&S(2\epsi_i) = \bigcup_{j=i+1}^{n}\{\epsi_i-\epsi_j,\epsi_i+\epsi_j\},~1\leq i\leq n,\\
&S(\epsi_i+\epsi_j)=\bigcup_{l=i+1}^{j-1}\{\epsi_i-\epsi_l,\epsi_l+\epsi_j\}
\cup\bigcup_{l=j+1}^n\{\epsi_i-\epsi_l,\epsi_j+\epsi_l\}\cup\\
&\bigcup_{l=j+1}^n\{\epsi_j+\epsi_l,\epsi_j-\epsi_l\}\cup S_{ij}
,~1\leq i<j\leq
n,\mbox{ where}\\
&S_{ij}=\left\{\begin{array}{ll}\{\epsi_i, \epsi_j\},&\mbox{if }\Phi=B_n,\\
%\{\epsi_i-\epsi_j, 2\epsi_j\}&\mbox{if }\Phi=C_n\\
\varnothing,&\mbox{if }\Phi=D_n.
\end{array}\right.
\end{array}\label{formula_sing_roots}
\end{equation}

\nota{Note that, given $\beta\in\Phi^+$, %$\beta\neq2\epsi_i$,
the column of $\beta$ contains exactly one of the root from each pair of $\beta$-singular roots whose sum is $\beta$. Precisely, if $\alpha, \beta, \gamma\in\Phi^+$,
$\beta=\alpha+\gamma$, then $\col(\beta)=\col(\alpha)$ if and only if $\alpha\succ\beta\succ\gamma$ and $\row(\gamma)=\row(\beta)$,
$\col(\gamma)=\row(\alpha)$, or $\col(\gamma)=-\row(\beta)$,
$\row(\gamma)=-\row(\alpha)$. Put
$$
S^+(\beta)=\{\alpha\in S(\beta)\mid\col(\alpha)=\col(\beta)\},~S^-(\beta)=\{\gamma\in S(\beta)\mid\col(\gamma)\neq\col(\beta)\},
$$
so that $S(\alpha)=S^+(\alpha)\cup
S^-(\alpha)$.\newpage

%If $\Phi=C_n$, $\beta=2\epsi_i$, and, say $\gamma=\epsi_i-\epsi_j\succ\alpha=\epsi_i+\epsi_j\succ\beta$ for some $j>i$, then $\col(\beta)=\col(\alpha)=\col(\gamma)=i$ and $\row(\gamma)=-\row(\alpha)=j$. In this case, we put
%$$S^+(\beta)=\bigcup_{j=i+1}^n\{\epsi_i+\epsi_j\},~S^-(\beta)=\bigcup_{j=i+1}^n\{\epsi_i-\epsi_j\},$$
%so that $S(\beta)=S^+(\beta)\cup S^-(\beta)$ again.

\label{nota_one_row_col}}

\defi{Let $D=\{\beta_1\succ\ldots\succ\beta_t\}\subset\Phi^+$ be an orthogonal rook placement
(more precisely, $\col(\beta_1)<\ldots<\col(\beta_t)$). We denote the four-elements set $\{\otimes,s, \bullet, + -\}$ by $\Sy$. The \emph{battleship}
is the map $\Sb=\Sb_D\colon\Phi^+\to\Sy$ defined by induction on $i, 0\leq i\leq t+1$, as follows.
\begin{enumerate}
\item If $\beta\in D$, then $\Sb(\beta)=\otimes$ (zero step).
\item Denote by $\Sy_j$
the set of all roots from $\Phi^+$ are mapped to minuses at the $j$th step (in particular, $\Sy_0=\varnothing$). Then at the $i$th step $i=1,\ldots,t$, we put $$\Sb(\gamma)=+,~\Sb(\delta)=-$$ for all $\gamma$, $\delta\in
S(\beta_i)$ such that $\gamma+\delta=\beta_i$,
$\gamma\in S^+(\beta_i)$ (so $\delta\in S^-(\beta_i)$) and $\gamma,
\delta\notin\cup_{j=0}^{i-1}\Sy_j$.
\item $\Sb(\alpha)=\bullet$ for all other roots (the $(t+1)$th step).
\end{enumerate}
The pair $(\Phi, \Sb)$ we will call the \emph{diagram},
associated with subsystem $D$.} \exam{Let $\Phi=B_6$,
$D=\{\epsi_1, \epsi_2+\epsi_5, \epsi_3-\epsi_6\}$. Below we schematically drew the battleship $\Sb_D$ as a filling of the diagram by symbols from $\Sy$.
\begin{center}\small
$\mymatrix{
\pho& \pho& \pho& \pho& \pho& \pho& \pho& \pho& \pho& \pho& \pho& \pho & \pho\\
\Top{2pt}\Rt{2pt}+& \pho& \pho& \pho& \pho& \pho& \pho& \pho& \pho& \pho& \pho& \pho & \pho\\
+& \Top{2pt}\Rt{2pt}+& \pho& \pho& \pho& \pho& \pho& \pho& \pho& \pho& \pho& \pho & \pho\\
+& +& \Top{2pt}\Rt{2pt}+& \pho& \pho& \pho& \pho& \pho& \pho& \pho& \pho& \pho & \pho\\
+& \bullet& \bullet& \Top{2pt}\Rt{2pt}\bullet& \pho& \pho& \pho& \pho& \pho& \pho& \pho& \pho & \pho\\
+& +& \otimes& -& \Top{2pt}\Rt{2pt}-& \pho& \pho& \pho& \pho& \pho& \pho& \pho & \pho\\
\otimes& -& -& -& -& \Top{2pt}\Rt{2pt}-& \pho& \pho& \pho& \pho& \pho& \pho & \pho\\
\bullet& +& \bullet& \bullet& -& \Top{2pt}\Lft{2pt}0& \Top{2pt}\Rt{2pt}\pho& \pho& \pho& \pho& \pho& \pho & \pho\\
\bullet& \otimes& -& -& \Top{2pt}\Lft{2pt}0& \pho& \pho& \Top{2pt}\Rt{2pt}\pho& \pho& \pho& \pho& \pho & \pho\\
\bullet& \bullet& \bullet& \Top{2pt}\Lft{2pt}0& \pho& \pho& \pho& \pho& \Top{2pt}\Rt{2pt}\pho& \pho& \pho& \pho & \pho\\
\bullet& \bullet& \Top{2pt}\Lft{2pt}0& \pho& \pho& \pho& \pho& \pho& \pho& \Top{2pt}\Rt{2pt}\pho& \pho& \pho & \pho\\
\bullet& \Top{2pt}\Lft{2pt}0& \pho& \pho& \pho& \pho& \pho& \pho& \pho& \pho& \Top{2pt}\Rt{2pt}\pho& \pho & \pho\\
\Top{2pt}0& \pho& \pho& \pho& \pho& \pho& \pho& \pho& \pho& \pho& \pho& \Top{2pt}\Rt{2pt}\pho & \pho\\
}$
\end{center} %For a step-by-step filling of the diagram, see Appendix.
\label{exam_sea_battle}}Next, we need the following simple Lemma. 
\lemmp{Let $\gamma\in\Phi^+$, $\Sb(\gamma)=+$. Then there exists exactly one root $\delta\in\Phi^+$ such that
$\gamma+\delta\in D$.\label{lemm_pro_plus}}{By definition of the battleship, there is a root $\delta\in\Phi^+$ such that $\Sb(\delta)=-$
and $\gamma+\delta=\beta_i$ for certain $\beta_i\in D$. Let us show that $\gamma+\alpha\notin D$ for all $\alpha\in\Phi^+$, $\alpha\neq\delta$.

Assume that there exist $\alpha$, $\beta_j\in\Phi^+$ such that
$\beta_j\in D$, $\gamma+\alpha=\beta_j$ and $\beta_j\neq\beta_i$. Since $\col(\gamma)=\col(\beta_i)$ by construction,
$\col(\gamma)\neq\col(\beta_j)$ (the orthogonal rook placement $D$ contains no more than one root from each column). Hence,
$\alpha\succ\beta_j\succ\gamma\succ\beta_i$, and
$\col(\alpha)=\col(\beta_j)$, $\row(\alpha)=\col(\gamma)$ (see
Remark~\ref{nota_one_row_col}: the situation
$\row(\alpha)=-\row(\gamma)$ is impossible, because then
$\col(\beta_i)=\col(\gamma)=-\row(\beta_j)$ and the roots
$\beta_i,\beta_j$ are not orthogonal). In particular, $j<i$.

The fact that the root $\gamma$ was not mapped to minus on the $j$th step of the battleship means that the root $\alpha$ was mapped to minus on the $l$th step for certain $l<j$. So, either
$\row(\alpha)=\row(\beta_l)$ or $\col(\alpha)=-\row(\beta_l)$.
However, in the first case $\row(\beta_l)=\col(\beta_i)$, and
$\beta_i$, $\beta_l$ are not orthogonal. In the second case,
$\row(\beta_l)=-\col(\beta_j)$, and $\beta_j$, $\beta_l$ are not orthogonal. Lemma is proved.}

It turns out that the battleship is a natural way to construct a polarization for the canonical form on an orbit associated with an orthogonal rook placement. More precisely, consider the set of roots which are not mapped to minuses via the battleship:
\begin{equation}
\Po=\Po_D=\{\alpha\in\Phi^+\mid\Sb(\alpha)\neq-\}.\label{formula_roots_pol}
\end{equation}
Denote by $\pt_D$ the subspace in $\ut$ spanned by the vectors
$e_{\alpha}$, $\alpha\in \Po$. \theop{Let $D\subset\Phi^+$ be an orthogonal rook placement\textup, $\xi$ is a tuple of non-zero scalars from~$\Fp_q$\textup,
$\Omega=\Omega_{D, \xi}\subset\ut^*$ be the corresponding associated coadjoint orbit\textup, and $f$ be the canonical form on it. Then $\pt=\pt_D$ is a polarization for $f$.\label{theo_pol}}{It consists of the three following steps.

1. The subspace $\pt$ is $f$-isotropic. It is enough to prove that the sum of any two roots from $\Po$ is not contained in $D$. But it follows immediately from the definition of the battleship.

2. Next, $\pt$ is a maximal $f$-isotropic subspace.
Indeed, assume that $\pt+\Fp_qy$ is $f$-isotropic for some $y\in\ut$. Let
$$y=\sum_{\delta\in\Phi^+\setminus \Po}y_{\delta}e_{\delta}\neq0,~y_{\delta}\in\Fp_q.$$ Since $y\neq0$, there exists
$\delta_0\in\Phi^+\setminus\Po$ such that $y_{\delta_0}\neq0$.
Since $\Sb(\delta_0)=-$, there exists $\gamma_0\in \Phi^+$ such that
$\Sb(\gamma_0)=+$ and $\gamma_0+\delta_0=\beta$ for certain
$\beta\in D$ (in particular, $\gamma_0\in\Po$, and so
$e_{\gamma_0}\in\pt$). It follows that $[e_{\gamma_0},
e_{\delta_0}]=ce_{\beta}$, $c\in\Fp_q^{\times}$. By Lemma \ref{lemm_pro_plus}, $\gamma_0+\delta\notin D$ for all
$\delta\in\Phi^+\setminus\{\delta_0\}$. Then
$$
f([e_{\gamma_0}, y])=f([e_{\gamma_0}, y_{\delta_0}e_{\delta_0}])=f(y_{\delta_0}ce_{\beta})=y_{\delta_0}c\xi_{\beta}\neq0,
$$
which contradicts the isotropy of the space $\pt+\Fp_qy$.

3. Finally, $\pt$ is a subalgebra. We postpone to Section~\ref{orbit_semidirect} the proof of this step %\ref{orbit_semidirect}
(see Proposition~\ref{prop_sub}), because it uses the semidirect decomposition of $U$.} %brbrbr there is another proof

\corop{Let $\Phi$\textup, $D$\textup, $\Omega$ be as in the theorem above. Then
$\dim\Omega$ is equal to the number of pluses and minuses in the diagram $(\Phi,
\Sb)$.\label{coro_dim_pol}}{The number of the pluses, by definition of the battleship, is equal to the number of minuses in $(\Phi, \Sb)$.
The dimension of an orbit is twice greater than the codimension of any polarization for any linear form on this orbit, but the construction of $\pt$ shows that $$\codim\pt=|(\Phi^+\setminus\Po)|
=|\{\alpha\in\Phi^+\mid\Sb(\alpha)=-\}|.$$ The proof is complete.} %Using this corollary, we will characterize the dimension of %\ref{formula_dimension}
%$\Omega$ in the 
%"inner"\ terms of Weyl group $W$ in Section 5 (see. Theorem
%\ref{theo_dim_orbit}). brbrbr Надо ли это?

In conclusion of this section, we note that the battleship first appeared in article \cite{IgnatevPanov09} under the name of ``admissible diagrams''. 
It allowed to describe all coadjoint orbits in the case $A_{n-1}$  when $n\leq6$, as well as particular series of orbits
(orbits of submaximal dimension, or orbits of depth 1, for example) for an arbitrary $n$. In the article 
\cite{Panov08}, the battleship was used for a complete description of orbits associated witch orthogonal rook placements for the case $\Phi=A_{n-1}$. In the paper \cite{IgnatevVenchakov24}, we used the battleship to classify so-called orbits of depth 2.

\bigskip
\sect{Semi-direct decomposition of the group $U$}
\label{orbit_semidirect}

In the previous section, we briefly described the connection between coadjoint
orbits and irreducible representations and characters of the group $U$. On the other hand, in the study of representations of finite groups,
an important role is played by the Mackey little group method based on the the notion of semi-direct
decomposition (see, for example, \cite{Lehrer}).

Let $\Gr$ be an arbitrary finite group, $A$, $B$ be its subgroups such that $\Gr$ is the semi-direct product of $A$ and $B$, i.e., $\Gr=AB$, $A\triangleleft\Gr$ and
$A\cap B=\{1\}$ (we denote $\Gr=A\rtimes B$). Assume additionally that $A$ is an abelian group and that $\psi$ is its irreducible character. The \emph{centralizer}
of $\psi$ in the group $B$ (or the \emph{little group}) is the subgroup of the form $$B^{\psi}=\{b\in
B\mid\psi\circ\tau_b=\psi\},$$ where $\tau_b\colon A\to A.~a\mapsto bab^{-1}.$

An arbitrary element $g\in\Gr$ can be uniquely represented as
$g=ab$, where $a\in A, b\in B$; this defines the maps
$$\pi^{\Gr}_{A}\colon\Gr\to A\colon g\mapsto a,~\pi^{\Gr}_{B}\colon\Gr\to B\colon g\mapsto b,~g\in\Gr.$$ Note that $\pi^{\Gr}_B$ is a group homomorphism, while $\pi^{\Gr}_A$, in general, is not. For an arbitrary subgroup $\wt B$ in~$B$ and arbitrary characters
$\psi\colon A\to\Cp$ and $\eta\colon\wt B\to\Cp$, one can define the character $\psi_0\eta_0$ of the group
$A\rtimes\wt B=A\wt B$, where
$$
\psi_0=\psi\circ\pi^{A\rtimes\wt
B}_{A},~\eta_0=\eta\circ\pi^{A\rtimes\wt B}_{\wt B}.
$$
Denote by $\Irr{\Gr}$ the set of irreducible characters of $\Gr$. The basic idea of the Mackey method, which allows us
to reduce the study of representations of the group $\Gr$ to the study of representations of
its subgroups $A$ and $B$, can be formulated as follows. \theo{Let $\Gr=A\rtimes B$
be a finite group with $A$ being an abelian group. Then each irreducible character $\chi$ of the group $\Gr$ has the form
\begin{equation}
\chi=\Ind{A\rtimes B^{\psi}}{G}{\psi_0\eta_0}
\label{formula_semi_direct}
\end{equation}
for certain $\psi\in\Irr{A}$\textup, $\eta\in\Irr{B}$. Conversely\textup, every function of the form \textup{(\ref{formula_semi_direct})} is an irreducible character of $\Gr$.}{\cite[Proposition
1.2]{Lehrer}}\label{theo_semi_direct}

Let us now show how to use these methods for our group $U$. %brbrbr не для C_n  и A_n
Let us denote by $\ut_1$ and $\vt$ the subalgebras of $\ut$ of the following form:
$$
\ut_1=\bigoplus_{\alpha\in
C_1}\Fp_qe_{\alpha},~\vt=\bigoplus_{\alpha\in\Phi^+\setminus
C_1}\Fp_qe_{\alpha}.
$$
Clearly, $\ut=\ut_1\oplus\vt$ as vector spaces and $\ut_1$ is a
 commutative ideal of $\ut$ (in this situation we will write $\ut=\ut_1\rtimes\vt$
and say that the algebra $\ut$ is the semi-direct product of its
subalgebras $\ut_1$ and $\vt$). Denote $U_1=\exp(\ut_1)$,
$V=\exp(\vt)$. \lemmp{The group $U$ is a semidirect product of the groups $U_1$ and $V$\textup, $U= U_1\rtimes V$\textup, and the subgroup
$U_1$ is abelian.\label{lemm_semidirect_U}}{It is well known that
%\cite[Corollary 2.5.17]{Carter}
$U=\prod_{\alpha\in\Phi^+}X_{\alpha}$, where positive roots
are taken in any fixed order and
$X_{\alpha}=\{\{x_{\alpha}(t),~t\in\Fp_q\}$ is the root subgroup corresponding to a positive root $\alpha\in\Phi^+$, where
$$
x_{\alpha}(t)=\exp(te_{\alpha})=\left\{\begin{array}{ll}
1_m+te_{\alpha},&\mbox{if }\row(\alpha)\neq0,\\
1_m+te_{\alpha}-\dfrac{t^2}{2}e_{-i, i},&\mbox{if
}\alpha=\epsi_i\in B_n.
\end{array}\right.
$$
It follows that $U=U_1V$, and
$$
U_1=\left\{\left(\begin{array}{cccccc}
1&0&0&\ldots&0&0\\
*&1&0&\ldots&0&0\\
*&0&1&\ldots&0&0\\
\vdots&\vdots&\vdots&\ddots&\vdots&\vdots\\
*&0&0&\ldots&1&0\\
*&*&*&\ldots&*&1\\
\end{array}\right)\in U\right\},~V=\left\{\left(\begin{array}{cccccc}
1&0&0&\ldots&0&0\\
0&1&0&\ldots&0&0\\
0&*&1&\ldots&0&0\\
\vdots&\vdots&\vdots&\ddots&\vdots&\vdots\\
0&*&*&\ldots&1&0\\
0&0&0&\ldots&0&1\\
\end{array}\right)\in U\right\}.
$$
The normality of the subgroup $U_1$ is now checked by simple matrix calculations. The commutativity of $U_1$ follows from the Backer--Campbell--Hausdorff formula (\ref{formula_Campbell_Hausdorf}) and the commutativity of $\ut_1$.} 
%\nota{We do not consider here the root system $C_n$ because this lemma is not true for it: the first is not a commutative subalgebra, because
%$(\epsi_1-\epsi_i)+(\epsi_1+\epsi_i)=2\epsi_1\in C_n^+$.}

Hence, we can use Theorem~\ref{theo_semi_direct} for $\Gr=U$, $A=U_1$, $B=V$. Now, let
$D=\{\beta_1\succ\ldots\succ\beta_t\}$ be a orthogonal rook placement in $\Phi^+$, $\Omega=\Omega_{D, \xi}\subset\ut^*$ be an associated coadjoint orbit, $f$ be the canonical form on it, and
$\chi=\chi_{\Omega}$
be the corresponding irreducible character of the group $U$. Our goal in the remainder of this section is to construct a decomposition of the form (\ref{formula_semi_direct}) for the character $\chi$.

Denote by $\psi=\psi_{D,\xi}$ the irreducible character of the group
$U_1$ of the form
\begin{equation}
\psi(x)=\theta(f(\ln(x)))=\left\{\begin{array}{ll}
\xi_{\beta_1}
x_{i,1},&\mbox{if }\col(\beta_1)=1,\row(\beta_1)=i,\\
1,&\mbox{if }\col(\beta_1)>1.
\end{array}\right.\label{formula_char_psi}
\end{equation}
(It is easy to see that, given $x=(x_{ij})\in U_1$, one has $\ln(x)=x-1_m-x_{-1, 1}e_{-1, 1}$. As usual, we denote the $(i,j)$th element of a matrix $x$ by $x_{ij}$.)

On the other hand, consider the following subgroup of $V$:
$$
\begin{array}{ll}
&V'=\prod_{\alpha\in \Phi_1^+}X_{\alpha},\mbox{ where}\\
&\Phi_1^+=\left\{\begin{array}{ll} \Phi^+\setminus C_1,&\mbox{if
}\col(\beta_1)>1,\\
\Phi^+\setminus(C_1\cup S^-(\beta_1)),&\mbox{if }\col(\beta_1)=1.
\end{array}\right.
\end{array}
$$
In particular, if $\col(\beta_1)>1$, then $V'=V$. It is clear that $V'=\exp(\vt')$, where $\vt'=\bigoplus_{\alpha\in \Phi_1^+}\Fp_qe_{\alpha}$.
Furthermore, if $\col(\beta_1)=1,\row(\beta_1)=i$, then
$$
V'=\{x\in V\mid x_{ij}=x_{-j,-i}=0,~1\prec j\prec i\}.
$$
\lemmp{Let $\col(\beta_1)=1$\textup, $\row(\beta_1)=i$. Then\textup, given $y\in\vt\setminus\vt'$\textup, there exists $z\in\ut_1$ such that $z_{i,1}\neq((\exp\ad{y})(z))_{i, 1}$.\label{lemm_z_i_1}}{Put
$$
y=\sum_{\alpha\in \Phi_1^+}y_{\alpha}e_{\alpha}+\sum_{\delta\in
S^-(\beta_1)}y_{\delta}e_{\delta}\in\vt\setminus\vt',~
y_{\alpha},~y_{\delta}\in \Fp_q.
$$
Let $\delta_0$ be the smallest root with respect to the order $\prec'$ (see
(\ref{formula_complete_orders})) from all $\delta\in S^-(\beta_1)$,
for which $y_{\delta}\neq0$.

Since $\Sb(\delta_0)=-$, by definition of the battleship,
there exists $\gamma_0\in\Phi^+$ such that $\Sb(\gamma_0)=+$ and
$\gamma_0+\delta_0=\beta_1$ (furthermore, in this case,
$\gamma_0\in C_1$, see Remark~\ref{nota_one_row_col}). Consider the element
$z=e_{\gamma_0}\in\ut_1$. By definition, $z_{i, 1}=0$. At the same time,
$$
(\exp\ad{y})(z)=z+\ad{y}z+\frac{1}{2}\ad{y}^2z+\ldots.
$$
Of course, for all $u\in\ut_1$ one has $u_{i, 1}=e_{\beta_1}^*(u)$,
hence, by Lemma~\ref{lemm_pro_plus}, $(\ad{y}z)_{i, 1}\neq0$, because\break $[e_{\gamma_0},e_{\delta_0}]=c_0e_{\beta_1}$ for some
$c_0\neq0$ and $[e_{\gamma_0},e_{\delta}]\notin \Fp_qe_{\beta_1}\setminus\{0\}$ for all
$\delta\in\Phi^+$, $\delta\neq\delta_0$.

Assume now that there exist $N>1$ and $\delta_1,~\ldots,~\delta_N\in \Phi^+\setminus C_1$ such that
$\gamma_0+\delta_1+\ldots+\delta_N=\beta_1$ and
$y_{\delta_i}\neq0$ for all $i=1,~\ldots,~N$. By \cite[Chapter VI, \S1,
Proposition 19]{Bourbaki03}, there exists $i_0$ such that
$\beta_1-\delta_{i_0}\in\Phi^+$; we may assume without loss of generality that $i_0=1$, i.e., $\delta_1\in S(\beta_1)$. But
$\delta_1\notin C_1$, hence $\delta_1\in S^-(\beta_1)$, because $S^+(\beta_1)\subset C_1$. So,
$\delta_1\succ'\delta_0$ by the choice of the root $\delta_0$. Thus, $\delta_0,~\delta_1\in S^-(\beta_1)$ and
$\delta_1\succ'\delta_0$; a simple case-by-case consideration shows that $\delta_1\nless\delta_0$ in the sense of the usual order on the roots
($\alpha\geq\beta$ if
$\alpha-\beta\in\langle\Phi^+\rangle_{\mathbb{Z}_{\gee0}}$). From the other side,
$$\delta_0=\beta_1-\gamma_0=\delta_1+\ldots+\delta_N>\delta_1,$$ which leads to a contradiction.
It follows that
$e_{\beta_1}^*(\ad{y}^N(z))=0$ for all $N>1$, so
$$((\exp\ad{y})(z))_{i,1}=c_0\neq0.$$ This means that $z\in\ut_1$ is an element we are looking for.} \lemmp{For an arbitrary orthogonal rook placement
$D\subset\Phi^+$ and an arbitrary $\xi$\textup, the subgroup $V'$
coincides with the centralizer of the character $\psi$ in the subgroup $V$\textup, i.e.\textup, $V'=V^{\psi}$.\label{lemm_stab_psi}}{It is clear if $\col(\beta_1)>1$, because in this case $V'=V^{\psi}=V$.
Let $$\col(\beta_1)=1, \row(\beta_1)=i.$$ Pick elements $x=\exp(y)\in V$ and $h=\exp(z)\in U_1$, where
$$
y=\sum_{\alpha\in\Phi^+\setminus
C_1}y_{\alpha}e_{\alpha}\in\vt,~z=\sum_{\gamma\in
C_1}z_{\gamma}e_{\gamma}\in\ut_1,~y_{\alpha},z_{\gamma}\in\Fp_q.
$$
Since $h=1_m+z+h_{-1,1}e_{-1,1}$, one has
$$
\begin{array}{ll}
xhx^{-1}&=\exp(y)(1_m+z+h_{-1,1}e_{-1,1})(\exp(y))^{-1}\\
&=1_m+\exp(y)z(\exp(y))^{-1}+h_{-1,1}\exp(y)e_{-1,1}(\exp(y))^{-1}.
\end{array}
$$
But $\exp(y)z(\exp(y))^{-1}=(\exp\ad{y})(z)$ and $e_{-1,1}$ is invariant under this conjugation. It follows that
$$
xhx^{-1}=1_m+h_{-1,1} e_{-1,1}+(\exp\ad{y})(z).
$$
In particular, by formula (\ref{formula_char_psi}), $$\psi(xhx^{-1})=\xi_{\beta_1}((\exp\ad{y})(z))_{i,
1}.$$ Since
$h_{i,1}=z_{i,1}$, Lemma~\ref{lemm_z_i_1} shows that $V'\supset
V^{\psi}$.

It remains to note that, for any $\alpha\in \Phi_1^+$, the subgroup
$X_{\alpha}$ centralizes the character $\psi$.
Hence, $V'\subset V^{\psi}$, which means that they coincide.}

So, we calculated the centralizer of the character $\psi$ in the subgroup $V$. Let us study its structure in more detail. Denote
$$
\Phi_2^+=\left\{\begin{array}{ll} \Phi_1^+\setminus
S^-(\epsi_1\pm\epsi_i),&\mbox{if
}\beta_1=\epsi_1\mp\epsi_i,~2\leq i\leq n,\\
\Phi_1^+&\mbox{otherwise},
\end{array}\right.
$$
and put
\begin{equation}
\vt_1=\bigoplus_{\alpha\in\Phi_1^+\setminus\Phi_2^+}\Fp_qe_{\alpha},~
\wt\ut=\bigoplus_{\alpha\in\Phi_2^+}\Fp_qe_{\alpha}.\label{formula_v_1}
\end{equation} It is easy to check that $\vt'=\vt_1\rtimes\wt\ut$, where $\vt_1$ is abelian. Arguing as in Lemma \ref{lemm_semidirect_U}, we see that
$$V'=V_1\rtimes\wt U$$ and $V_1$ is abelian, where
$$V_1=\exp(\vt_1)=\prod_{\alpha\in\Phi_1^+\setminus\Phi_2^+}X_{\alpha},~\wt U=\exp(\wt\ut)=\prod_{\alpha\in\Phi_2^+}X_{\alpha}$$ (if
$\Phi_1^+=\Phi_2^+$ then $V_1=\{1_m\}$ and $\wt U=V'$).

Moreover, we notice that $\wt U$ is isomorphic to the maximal unipotent subgroup of the orthogonal group of smaller rank. More precisely, set
$$
\wt\Phi=\left\{\begin{array}{ll}B_{n-1},&\mbox{if }\Phi=B_n,\col(\beta_1)>1,\\
D_{n-1},&\mbox{if }\Phi=D_n,~\col(\beta_1)>1,\\
D_{n-2},&\mbox{if }\Phi=D_n,~\col(\beta_1)=1,\\
B_{n-2},&\mbox{if }\Phi=B_n,~\col(\beta_1)=1,~\row(\beta_1)\neq0,\\
D_{n-1},&\mbox{if }\Phi=B_n,~\col(\beta_1)=1,~\row(\beta_1)=0.
\end{array}\right.
$$
\lemmp{There is an isomorphism
$\wt\ut\cong\ut(\wt\Phi)$.\label{lemm_isomorph_RS}}{It is enough to construct a one-to-one map
$\pi\colon\Phi_2^+\to\wt\Phi^+$, which can be continued to an isometric isomorphism
$\langle\Phi_2^+\rangle_{\Rp}\to\langle\wt\Phi^+\rangle_{\Rp}$.
Note that $\Phi_2^+$ is obtained from $\Phi^+$ by removing some rows and columns. Namely,
$$
\Phi_2^+=\left\{\begin{array}{ll}\Phi_1^+=\Phi^+\setminus C_1,&\mbox{if }\col(\beta_1)>1,\\
\Phi^+\setminus(C_1\cup R_0),&\mbox{if
}\col(\beta_1)=1,~\row(\beta_1)=0,\\
\Phi^+\setminus(C_1\cup R_i\cup R_{-i}\cup C_i),&\mbox{if
}\col(\beta_1)=1,~\row(\beta_1)=i\neq0.
\end{array}\right.
$$
Define the number $\wt m$ for $\wt\Phi$ similarly to the number $m$ for
$\Phi$ (see Section 2%\ref{main_definitions}
). If it is even then
put $\wt n=\wt m/2$, otherwise put $\wt n=(\wt m+1)/2$. Let us numerate columns of roots in $\Phi_2^+$ from $1$ to $\wt n$, and their rows from
$-\wt n$ to $\wt n$ (except $0$ if $\wt m$ is even). We obtain the required map $\pi$. }

We will denote the corresponding isomorphism
$\wt\ut\to\ut(\wt\Phi)$ defined by the rule $e_{\alpha}\mapsto
e_{\pi(\alpha)}$, $\alpha\in\Phi_2^+$, by the same letter $\pi$.

Now we are ready to prove that the space $\pt$ constructed via the battleship is a subalgebra. This fact announced in Section~\ref{sect:sea_battle_pol} proves that $\pt$ is a polarization for the canonical form on an orbit associated with an orthogonal rook placement, see
Theorem~\ref{theo_pol}. \propp{Let $\Phi=B_n$ or $D_n$\textup,
$D\subset\Phi^+$
be an arbitrary orthogonal rook placement\textup, and $\pt$ be the subspace constructed by the battleship $\Sb_D$. Then
$\pt$ is a subalgebra in $\ut$.\label{prop_sub}}{We will use the induction on $n$. The base of induction can be checked directly. Denote by $\wt D$ the orthogonal rook placement $\pi(D\setminus
C_1)\subset\wt\Phi^+$. Let $\wt\Sb$ be the corresponding battleship on $\wt\Phi^+$. Then
$\Sb{\mid}_{\Phi_2^+}=\wt\Sb\circ\pi$, hence $\wt\pt=\pt\cap\wt\ut$
is a subalgebra in $\wt\ut$ (and so in $\ut$) by the inductive assumption as the preimage $\pi^{-1}(\pt_{\wt D})$ of the subalgebra
$\pt_{\wt D}\subset\ut(\wt\Phi)$.

It follows from the definition of the battleship (see Section~\ref{sect:sea_battle_pol}) that if a root $\alpha$ belongs to the first column of $\Phi^+$ then
$\Sb(\alpha)\neq-$, so $\ut_1\subset\pt$. Denote
$\pt_1=\pt\cap\vt_1=\bigoplus_{\alpha\in\Phi_1^+\setminus\Phi_2^+}e_{\alpha}$.
Obviously, $\pt=\ut_1\oplus\pt_1\oplus\wt\pt$ as vectors
spaces. Since $\ut_1\triangleleft\ut$, $\pt_1$ is abelian and $\wt\pt$ is a subalgebra of $\ut$, it suffices to show that $[\wt\pt,~\pt_1]\subset\pt$ to complete the proof.
Consider three possible cases for $\beta_1$.

1. If $\col(\beta_1)>1$ or $\Phi=B_n$ and $\beta_1=\epsi_1$ (i.e., 
$\col(\beta_1)=1$, $\row(\beta_1)=0$), then
$\Phi_1^+\setminus\Phi_2^+=\varnothing$ and, consequently, $\pt_1=0$.

2. Let $\col(\beta_1)=1$, $\row(\beta_1)=-i$ (i.e,
$\beta_1=\epsi_1+\epsi_i$). Then
$$
\Phi_2^+\setminus\Phi_1^+=S^-(\epsi_1-\epsi_i)=\{\epsi_j-\epsi_i,~2\leq
j\leq i-1\}.
$$
Assume that $\alpha\in S^-(\epsi_1-\epsi_i)$, $\delta\in
\Po\cap\Phi_2^+$ and $\alpha+\delta=\gamma\in\Phi^+$. (Recall that
$\pt=\bigoplus_{\alpha\in\Po}\Fp_qe_{\alpha}$, see the definition of $\Po$ in (\ref{formula_roots_pol}).) Since $\Phi_2^+$ (and even $\Phi_1^+$)
does not contain $C_1$, $R_{-i}$ and $C_i$, and $\alpha$  has the form $\epsi_j-\epsi_i$ for a certain $2\leq j\leq i-1$, the root $\gamma$ equals $\epsi_l-\epsi_i$ for a certain $2\leq l\leq i-2$. This means that $\gamma\in\Po$, and so $e_{\gamma}\in\pt$.

3. Finally, let $\col(\beta_1)=1$, $\row(\beta_1)=i$ (i.e.,
$\beta_1=\epsi_1-\epsi_i$). Here
$$
\Phi_2^+\setminus\Phi_1^+=S^-(\epsi_1+\epsi_i)=\{\epsi_j+\epsi_i,~2\leq
j\leq i-1\}\cup C_i.
$$
The proof is rather similar to the previous case. Precisely, let
$\alpha\in S^-(\epsi_1+\epsi_i)$, $\delta\in\Po\cap\Phi_2^+$ and
$\alpha+\delta=\gamma\in\Phi^+$. Since $\Phi_2^+$ (and even
$\Phi_1^+$) does not contain $R_i$, so $\epsi_i$ enters the root $\gamma$ with positive coefficient. Hence, $\gamma\in R_{-i}\cup C_i\subset\Po$. Thus, $e_{\gamma}\in\pt$. The proof is complete.}

Now, finally, we can present a decomposition of the form
(\ref{formula_semi_direct}) for the irreducible character
associated with an orthogonal rook placement. Recall that we have constructed
semi-direct decomposition of $U=U_1\rtimes V$, defined the irreducible
character $\psi$ of the group $U_1$ (cf. (\ref{formula_char_psi})) and found
its centralizer $V^{\psi}=V'=V_1\rtimes\wt U$. In this case, $\wt\ut=\mathrm{Lie}(\wt U)\cong\ut(\wt\Phi)$, so
the irreducible characters of the group $\wt U$ are in
one-to-one correspondence with the coadjoint orbits in $\wt\ut^*$ (cf. (\ref{formula_char_orbit})).

As above, let  $D=\{\beta_1\succ\ldots\succ\beta_t\}$~ be
an orthogonal rook placement in $\Phi^+$, $\xi=(\xi_{\beta})_{\beta\in D}$,
$\Omega=\Omega_{D,\xi}$ be the associated orbit, and $\chi$
be the corresponding irreducible character of the group $U$. As in the
proof of Proposition \ref{prop_sub}, we denote $\wt D
=\pi(D\setminus C_1)\subset\wt\Phi^+$,
$\wt\Omega=\Omega_{\wt D,\wt\xi}$, and let $\wt\chi$ be the corresponding irreducible character of the group $\wt U$. Here
$$
\wt\xi=\left\{\begin{array}{ll} \xi\in(\Fp_q^{\times})^t,&\mbox{if }\col(\beta_1)>1,\\
(\xi_{\beta_2},\ldots,\xi_{\beta_t})\in(\Fp_q^{\times})^{t-1},&\mbox{if
}\col(\beta_1)=1.
\end{array}\right.
$$

\theop{Let $\pi^{U_1\rtimes V'}_{\wt U}=\pi^{V'}_{\wt
U}\circ\pi^{U_1\rtimes V'}_{V'}$. Then \begin{equation}
\chi=\Ind{U_1\rtimes V'}{U}{((\psi\circ\pi^{U_1\rtimes
V'}_{U_1})\cdot(\wt\chi\circ\pi^{U_1\rtimes V'}_{\wt
U}))}.\label{formula_char_decompose}
\end{equation}.\label{theo_char_decompose}}{Put $P=\exp(\pt)\subset U, \wt
P=\exp(\wt\pt)\subset\wt U$, then $P=U_1\rtimes(V_1\rtimes\wt
P)$. Let us denote the function on the right-hand side of
(\ref{formula_char_decompose}) by $\eta$. By Theorem
\ref{theo_char_orbit},
$$
\chi=\Ind{P}{U}{(\theta\circ f\circ\ln)},\quad \wt\chi=\Ind{\wt
P}{\wt U}{(\theta\circ f\circ\ln)}.
$$
We will break down further reasoning into a series of steps.

1. We claim that if $y\in P$ then $f(\ln\pi^{P}_{V_1\rtimes\wt
P}(y))=f(\ln\pi^P_{\wt P}(y))$. Indeed, it is evident that
$$
\pi^{P}_{V_1\rtimes\wt P}(y)=(\pi^{V_1\rtimes\wt
P}_{V_1}\circ\pi^{P}_{V_1\rtimes\wt P})(y)\cdot(\pi^{V_1\rtimes\wt
P}_{\wt P}\circ\pi^{P}_{V_1\rtimes\wt
P})(y)=\pi^P_{V_1}(y)\cdot\pi^P_{\wt P}(y).
$$
Let $u=\ln\pi^P_{V_1}(y)\in\pt$, $v=\ln\pi^P_{\wt P}(y)\in\pt$, then, from (\ref{formula_Campbell_Hausdorf}), one has
$$\pi^{P}_{V_1\rtimes\wt P}(y)=\exp(u)\exp(v)=\exp(u+v+\tau),
\quad\mbox{where }\tau\in[\pt,\pt],$$ so
$\ln\pi^{P}_{V_1\rtimes\wt P}(y)=\ln\pi^P_{V_1}(y)+\ln\pi^P_{\wt
P}(y)+\tau$. But $\pt$ is a polarization for $f$, so $f(\tau)=0$. At the same time, $\ln\pi^P_{V_1}(y)\in\vt_1$ and, clearly,
$f{\mid}_{\vt_1}\equiv0$ (see the definition of $\vt_1$ in
(\ref{formula_v_1})). It follows that
$$
f(\ln\pi^{P}_{V_1\rtimes\wt
P}(y))=f(\ln\pi^P_{V_1}(y))+f(\ln\pi^P_{\wt
P}(y))+f(\tau)=f(\ln\pi^P_{\wt P}(y)).
$$

2. We claim that if $y\in P$, then $f(\ln y-\ln\pi^P_{\wt
P}(y)-\ln\pi^P_{U_1}(y))=0$. The proof is similar to the previous step. Indeed, if $u=\ln\pi^P_{V_1\rtimes \wt P}(y)\in\pt$,
$v=\ln\pi^P_{U_1}(y)\in\pt$ then, according to
(\ref{formula_Campbell_Hausdorf}),
$$
y=\pi^P_{U_1}(y)\cdot\pi^P_{V_1\rtimes\wt
P}(y)=\exp(u)\exp(v)=\exp(u+v+\tau),\quad\mbox{where }\tau\in[\pt,\pt].
$$
Hence, $\ln y=\ln\pi^P_{U_1}(y)+\ln\pi^P_{V_1\rtimes\wt P}+\tau$.
But $\pt$ is a polarization for $f$, so $f(\tau)=0$ and, consequently,
$f(\ln\pi^{P}_{V_1\rtimes\wt P}(y))=f(\ln\pi^P_{\wt P}(y))$.

3. The next claim is as follows. Let $\Gr= A\rtimes B$ be a finite group, $C$ be a subgroup in $B$ and $\lambda$ be a complex representation of the group $C$. Then
$\Ind{A\rtimes C}{G}{(\lambda\circ\pi^{A\rtimes
C}_C)}=(\Ind{C}{B}{\lambda})\circ\pi^{\Gr}_B$. This is proved in \cite[Proposition
1.2]{Lehrer}.

4. Let us apply step 3 to the case $\Gr=V_1\rtimes\wt U$, $A=V_1$, $B=\wt
U$, $C=\wt P$, $\lambda=\theta\circ f\circ\ln$. We obtain that
$$
\Ind{\wt P}{\wt U}{(\theta\circ f\circ\ln)}\circ\pi^{V_1\rtimes\wt
U}_{\wt U}=\Ind{V_1\rtimes\wt P}{V_1\rtimes\wt U}{(\theta\circ
f\circ\ln\circ\pi^{V_1\rtimes\wt P}_{\wt P})}.
$$
Now we apply step 3 to the case of $\Gr=U_1\rtimes V'$,
$A=U_1$, $B=V'=V_1\rtimes\wt U$, $C=V_1\rtimes\wt P$,
$\lambda=\theta\circ f\circ\ln\circ\pi^{V_1\rtimes\wt P}_{\wt P}$.
We have
$$
\begin{array}{ll}
\Ind{V_1\rtimes\wt P}{V_1\rtimes\wt U}{(\theta\circ
f\circ\ln\circ\pi^{V_1\rtimes\wt P}_{\wt P})}\circ\pi^{U_1\rtimes
V'}_{V_1\rtimes\wt U}&=\Ind{U_1\rtimes(V_1\rtimes\wt P)}{U_1\rtimes
V'}{(\theta\circ f\circ\ln\circ\pi^{V_1\rtimes\wt P}_{\wt
P}\circ\pi^{U_1\rtimes (V_1\rtimes \wt P)}_{V_1\rtimes\wt P})}=\\
&=\Ind{P}{U_1\rtimes V'}{(\theta\circ f\circ\ln\circ\pi^P_{\wt P})}.
\end{array}
$$
Hence, keeping in mind that $V'=V_1\rtimes\wt U$, we obtain
$$
\begin{array}{ll}
\Ind{\wt P}{\wt U}{(\theta\circ f\circ\ln)}\circ\pi^{U_1\rtimes
V'}_{\wt U}&=(\Ind{\wt P}{\wt U}{(\theta\circ
f\circ\ln)})\circ\pi^{V_1\rtimes\wt U}_{\wt U}\circ\pi^{U_1\rtimes
V'}_{V_1\rtimes \wt U}=\\
&=\Ind{P}{U_1\rtimes V'}{(\theta\circ f\circ\ln\circ\pi^P_{\wt P})}.
\end{array}
$$

5. One can rewrite the right-hand side of (\ref{formula_char_decompose}) taking into account
step 4 in the following way:
$$
\begin{array}{ll}
\eta&=\Ind{U_1\rtimes V'}{U}{((\psi\circ\pi^{U_1\rtimes
V'}_{U_1})\cdot(\wt\chi\circ\pi^{U_1\rtimes V'}_{\wt
U}))}=\\
&=\Ind{U_1\rtimes V'}{U}{((\theta\circ f\circ\ln\circ\pi^{U_1\rtimes
V'}_{U_1})\cdot(\Ind{\wt P}{\wt U}{(\theta\circ
f\circ\ln)}\circ\pi^{U_1\rtimes V'}_{\wt U}))}=\\
&=\Ind{U_1\rtimes V'}{U}{((\theta\circ f\circ\ln\circ\pi^{U_1\rtimes
V'}_{U_1})\cdot\Ind{P}{U_1\rtimes V'}{(\theta\circ
f\circ\ln\circ\pi^P_{\wt P})})}.
\end{array}
$$

6. We claim that the following equality holds:
$$
\Ind{P}{U_1\rtimes V'}{(\theta\circ f\circ\ln)}=(\theta\circ
f\circ\ln\circ\pi^{U_1\rtimes V'}_{U_1})\cdot\Ind{P}{U_1\rtimes
V'}{(\theta\circ f\circ\ln\circ\pi^P_{\wt P})}.
$$
Indeed, let $H\subset U_1\rtimes V'$ be an arbitrary complete
representative system of $(U_1\rtimes V')/P$ (in fact, one can choose $H$ to be a subset of $V'$, because $U_1\subset P$). Then, for any $x\in U_1\rtimes V'$,
$$
\begin{array}{ll}
&((\theta\circ f\circ\ln\circ\pi^{U_1\rtimes
V'}_{U_1})\cdot\Ind{P}{U_1\rtimes V'}{(\theta\circ
f\circ\ln\circ\pi^P_{\wt P})})(x)=\\
&=\theta(f(\ln\pi^{U_1\rtimes V'}_{U_1}(x)))\cdot\sum_{h\in H\mid
h^{-1}xh\in
P}\theta(f(\ln\pi^{P}_{\wt P}(h^{-1}xh)))=\\
&=\sum_{h\in H\mid h^{-1}xh\in P}\theta(f(\ln\pi^{U_1\rtimes
V'}_{U_1}(x)+\ln\pi^P_{\wt P}(h^{-1}xh))).
\end{array}
$$
But $h\in B$ implies that $h^{-1}\pi^{\Gr}_{A}(x)h=\pi^{\Gr}_A(h^{-1}xh)$
for an arbitrary finite group $\Gr=A\rtimes B$ and arbitrary
$x\in G$, $h\in B$. In our situation, this means that
$$\psi(X)=f(\ln\pi^{U_1\rtimes V'}_{U_1}(x))=f(\ln(h^{-1}\pi^{U_1\rtimes
V'}_{U_1}(x)h))=f(\ln\pi^{U_1\rtimes V'}_{U_1}(h^{-1}xh)),$$
because $h\in H\subset V' = V^{\psi}$, see Lemma
\ref{lemm_stab_psi}.

In addition, for $h^{-1}xh\in P$, we have $\pi^{U_1\rtimes
V'}_{U_1}(h^{-1}xh)=\pi^{P}_{U_1}(h^{-1}xh)$. Using step 2,
we obtain
$$
\begin{array}{ll}
&((\theta\circ f\circ\ln\circ\pi^{U_1\rtimes
V'}_{U_1})\cdot\Ind{P}{U_1\rtimes V'}{(\theta\circ
f\circ\ln\circ\pi^P_{\wt P})})(x)=\\
&=\sum_{h\in H\mid h^{-1}xh\in
P}\theta(f(\ln\pi^{P}_{U_1}(x)+\ln\pi^P_{\wt P}(h^{-1}xh)))=\\
&=\sum_{h\in H\mid h^{-1}xh\in P}\theta(f(\ln
h^{-1}xh))=\Ind{P}{U_1\rtimes V'}{(\theta\circ f\circ\ln)}.
\end{array}
$$

7. Combining steps 5 and 6, we conclude that $$\eta=\Ind{U_1\rtimes
V'}{U}{\Ind{P}{U_1\rtimes V'}{(\theta\circ
f\circ\ln)}}=\Ind{P}{U}{(\theta\circ f\circ\ln)}=\chi,$$ which completes the proof.}

So, the computation of the irreducible character associated with
an arbitrary orthogonal rook placement can be reduced, in fact, to the
computation of an irreducible character in the maximal unipotent subgroup
of an orthogonal group of smaller rank. This inductive procedure
will allow us to prove the formula for the dimension of an orbit,
associated with an orthogonal rook placement.

\sect{Formula for the dimension}\label{formula_dimension}

We recall all the notation from the previous sections. In particular,
$D=\{\beta_1\succ\ldots\succ\beta_t\}$ is an orthogonal rook placement in
$\Phi^+$, $\Omega=\Omega_{D,\xi}$ is an associated coadjoint orbit of the group $U$, $\chi$ is the corresponding irreducible character. Let $\sigma$ be the involution in the Weyl group $W(\Phi)$ defined by $\Supp{\sigma}=D$ (see (\ref{formula_ortog_decompose})). Our goal in
this section is to characterize the dimension of $\Omega$ in the ``interior" terms of the Weyl group $W$ (more precisely, in terms of the involution $\sigma$), combining
Theorem \ref{theo_char_decompose} and Corollary \ref{coro_dim_pol}.

Recall the definition of the root system $\wt\Phi$. Let $\wt W$
be its Weyl group, let $\wt\sigma\in\wt W$ be the involution with the support $$\Supp{\wt\sigma}=\wt D=\pi(D\setminus C_1)$$ (the map
$\pi\colon\Phi_2^{+}\to\wt\Phi^{+}$ is defined in Lemma
\ref{lemm_isomorph_RS}). Let $\sigma_2\in W$ be the involution with the support $D\cap\Phi_2^+$. For an arbitrary involution $\tau\in
W(\Phi)$, put
$$
\Phi_{\tau}=\{\alpha\in\Phi^+\mid\tau(\alpha)<0\}.
$$
The calculation of $\dim\Omega$ is essentially based on the comparison of the sets
$\Phi_{\sigma}$, $\Phi_{\sigma_2}$ and $\wt\Phi_{\wt\sigma}$.

We will need a few more notations. For an
involution $\tau\in W(\Phi)$, whose support is an orthogonal rook placement, we denote by $l(\tau)$ (respectively, by $s(\tau)$)
the length of the reduced (i.e. the shortest) decomposition of $\tau$ in the product of simple (respectively, of arbitrary) reflections. It is well known that
$l(\tau)=|\Phi_{\tau}|$. Clearly, $s(\tau)=|\Supp{\tau}|$. Let us also determine the number $d(\tau)$ by the formula
$$
d(\tau)=\left\{\begin{array}{ll}
\#\{\beta\in\Supp{\tau}\mid\col(\beta)>i,\row(\beta)<0\},&\mbox{if
}\Phi=B_n,~\epsi_i\in\Supp{\tau},\\
0,&\mbox{else.}
\end{array}\right.
$$
Since $\Supp{\tau}$ is an orthogonal rook placement (cf.
(\ref{formula_basic_supp})), it can not contain more than one
root of the form $\epsi_i$, because such roots lie in $R_0$
so the number $d(\tau)$ is well-defined. \exam{Let
$\Phi=B_6$,
$D=\Supp{\sigma}=\{\epsi_1,\epsi_2+\epsi_6,\epsi_3+\epsi_5\}$. We draw the diagram $(\Phi,\Sb)$ corresponding to the battleship
$\Sb=\Sb_D$ below. Here $d(\sigma)=2$.
\begin{center}\small
$\mymatrix{
\pho& \pho& \pho& \pho& \pho& \pho& \pho& \pho& \pho& \pho& \pho& \pho & \pho\\
\Top{2pt}\Rt{2pt}+& \pho& \pho& \pho& \pho& \pho& \pho& \pho& \pho& \pho& \pho& \pho & \pho\\
+& \Top{2pt}\Rt{2pt}+& \pho& \pho& \pho& \pho& \pho& \pho& \pho& \pho& \pho& \pho & \pho\\
+& +& \Top{2pt}\Rt{2pt}+& \pho& \pho& \pho& \pho& \pho& \pho& \pho& \pho& \pho & \pho\\
+& +& \bullet& \Top{2pt}\Rt{2pt}\bullet& \pho& \pho& \pho& \pho& \pho& \pho& \pho& \pho & \pho\\
+& \bullet& \bullet& \bullet& \Top{2pt}\Rt{2pt}\bullet& \pho& \pho& \pho& \pho& \pho& \pho& \pho & \pho\\
\otimes& -& -& -& -& \Top{2pt}\Rt{2pt}-& \pho& \pho& \pho& \pho& \pho& \pho & \pho\\
\bullet& \otimes& -& -& -& \Top{2pt}\Lft{2pt}0& \Top{2pt}\Rt{2pt}\pho& \pho& \pho& \pho& \pho& \pho & \pho\\
\bullet& \bullet& \otimes& -& \Top{2pt}\Lft{2pt}0& \pho& \pho& \Top{2pt}\Rt{2pt}\pho& \pho& \pho& \pho& \pho & \pho\\
\bullet& \bullet& \bullet& \Top{2pt}\Lft{2pt}0& \pho& \pho& \pho& \pho& \Top{2pt}\Rt{2pt}\pho& \pho& \pho& \pho & \pho\\
\bullet& \bullet& \Top{2pt}\Lft{2pt}0& \pho& \pho& \pho& \pho& \pho& \pho& \Top{2pt}\Rt{2pt}\pho& \pho& \pho & \pho\\
\bullet& \Top{2pt}\Lft{2pt}0& \pho& \pho& \pho& \pho& \pho& \pho& \pho& \pho& \Top{2pt}\Rt{2pt}\pho& \pho & \pho\\
\Top{2pt}0& \pho& \pho& \pho& \pho& \pho& \pho& \pho& \pho& \pho& \pho& \Top{2pt}\Rt{2pt}\pho & \pho\\
}$
\end{center}\label{exam_inv_B_6}}

To prove the main result of this section, we will examine $l(\sigma)$ and $l(\wt\sigma)$
for different cases of the root~$\beta_1$. It is obvious that
if $\col(\beta_1)>1$ then $l(\sigma)=l(\wt\sigma)$, because
$\Phi_2^+=\Phi_1^+=\Phi^+\setminus C_1$. Let $\col(\beta_1)=1$, $\row(\beta_1)=\pm i\neq0$ (i.e.,
$\beta_1=\epsi_1\mp\epsi_i$). Since
$\pi^{-1}(\wt\Phi^+)=\Phi_2^+$ and $\Phi_2^+$ does not contain $C_1, C_i,
R_i$ and $R_{-i}$, for $\alpha\in\Phi_2^+$, the conditions
$\sigma(\alpha)<0$, $\sigma_2(\alpha)<0$ and
$\wt\sigma(\pi(\alpha))<0$ are equivalent. It is also clear that if $\alpha\in C_1$ then $\sigma_2(\alpha)>0$,
hence $\sigma(\alpha)<0$ if and only if $\sigma_2(\alpha)\in
S^+(\beta_1)\cup\{\beta_1\}$, i.e., $\sigma_2(\alpha)\preceq'\beta_1$ (or, equivalently, $\row(\sigma_2(\alpha))\preceq\row(\beta_1)$). Similarly, if
$\alpha=\epsi_i\in R_0$, then $\sigma_2(\alpha)=\alpha>0$, and
$\sigma(\alpha)>0$ if and only if $\row(\beta_1)>0$. This means that $\Phi_{\sigma_2}$ can be obtained from
$\pi^{-1}(\wt\Phi_{\wt\sigma})$ by adding the roots from $(R_i\cup
R_{-i}\cup C_i)\setminus (C_1\cup R_0)$ which become negative under the action of the $\sigma_2$ (we put
$R_0=\varnothing$ in the case when $\Phi=D_n$).\newpage \lemmp{Let
$\sigma,\wt\sigma$ be as above\textup, and
$\beta_1=\epsi_1-\epsi_i$\textup, $2\leq i\leq n$. Then
$l(\sigma)=l(\wt\sigma)+2i-3$.\label{lemm_sigma_i_m_j}}{According to the previous remarks, $\Phi_{\sigma}\cap C_1$ consists from $(i-2)$ roots of the form $\sigma_2^{-1}(\epsi_1-\epsi_j)$, $1<j<i$ (they belong to $S^+(\beta_1)$), together with the root $\beta_1$. Assume that $\alpha=\epsi_j+\epsi_i\in R_{-i}$, $1<j<i$. The number $i$ does not appear among the rows and columns of the roots $\beta_2,\ldots,\beta_t$, hence
$\sigma_2(\alpha)=\epsi_i\pm\epsi_l$ for some $l$. In any
case, $\sigma(\alpha)=\epsi_1\pm\epsi_l>0$.

Assume now, that $\alpha=\epsi_i\pm\epsi_j\in C_i$, $i<j\leq
n$. Then $\sigma_2(\alpha)=\epsi_i\pm\epsi_l$ for some $l$
(signs are independent), and $\sigma(\alpha)=\epsi_1\pm\epsi_l>0$. Finally,
for $\alpha=\epsi_j-\epsi_i\in R_i$, $1<j<i$, we have
$\sigma(\alpha)=-\epsi_1\pm\ldots<0$ anyway.

Hence, all roots from $R_{-i}\cup C_i$ are mapped by the involution
$\sigma$ to positive roots, and all roots from $R_i$ are mapped to negative. Thus,
$$
\begin{array}{ll}
l(\sigma)&=|\Phi_{\sigma}|=|\pi^{-1}(\wt\Phi_{\wt\sigma})\cup
R_i\cup
S^+(\beta_1)|=|\wt\Phi_{\wt\sigma}|+|R_i|+|S^+(\beta_1)|=\\
&=l(\wt\sigma)+(i-1)+(i-2)=l(\wt\sigma)+2i-3,
\end{array}
$$
as required. } Recall the definition of $m$ from Section \ref{main_definitions}.\lemmp{Let $\sigma,\wt\sigma$
be as above\textup, and $\beta_1=\epsi_1+\epsi_i$\textup, $2\leq i\leq
n$. Then
$l(\sigma)=l(\wt\sigma)+2(m-i)-3$.\label{lemm_sigma_i_p_j}}{According to the previous remarks, $\Phi_{\sigma}\cap C_1$ consists of 
$(m-i-2)$ roots of the form $\sigma_2^{-1}(\epsi_1\pm\epsi_j)$, $1<j<i$ (they belong to $S^+(\beta_1)$; in the case of $\Phi=B_n$, the root $\epsi_1$ should also be counted), together with the root $\beta_1$. Assume that $\alpha=\epsi_j-\epsi_i\in R_{-i}$, $1<j<i$. The number $i$ does not appear among the rows and columns of the roots $\beta_2,\ldots,\beta_t$, hence
$\sigma_2(\alpha)=-\epsi_i\pm\epsi_l$ for some $l$. In any
case, $\sigma(\alpha)=\epsi_1\pm\epsi_l>0$.

Assume now that $\alpha=\epsi_i\pm\epsi_j\in C_i$, $i<j\leq
n$. Then $\sigma_2(\alpha)=\epsi_i\pm\epsi_l$ for some $l$
(signs are independent), and $\sigma(\alpha)=-\epsi_1\pm\epsi_l<0$.
Finally, when $\alpha=\epsi_j+\epsi_i\in R_{-i}$, $1<j<i$, we have
$\sigma(\alpha)=-\epsi_1\pm\ldots<0$ anyway.

Hence, all roots from $R_i$ are mapped by the involution
$\sigma$ to positive roots, and all roots from $R_{-i}\cup C_i$ are mapped to negative once. Thus,
$$
\begin{array}{ll}
l(\sigma)&=|\Phi_{\sigma}|=|\pi^{-1}(\wt\Phi_{\wt\sigma})\cup
R_{-i}\cup C_i\cup
S^+(\beta_1)|=|\wt\Phi_{\wt\sigma}|+|R_{-i}|+|C_i|+|S^+(\beta_1)|=\\
&=l(\wt\sigma)+(i-1)+(m-2i)+(m-i-2)=l(\wt\sigma)+2(m-i)-3,
\end{array}
$$
as required.}

Let us now consider the case when $\beta_1=\epsi_1$ (i.e.,
$\row(\beta_1)=0$). \lemmp{Let $\Phi=B_n$\textup, $\sigma,\wt\sigma$ be as above\textup, and $\beta_1=\epsi_1$. Then
$l(\sigma)=l(\wt\sigma)+2(n+d(\sigma))-1$.\label{lemm_sigma_i}}{Here
$\Phi_2^+=\Phi_1^+=\Phi^+\setminus(C_1\cup R_0)$. The reflection
$r_{\beta_1}$ acts on $\Phi_2^+$ identically, therefore, for
$\alpha\in\Phi_2^+$, the conditions $\sigma(\alpha)<0$, $\sigma_2(\alpha)<0$
and $\wt\sigma(\pi(\alpha))<0$ are equivalent. It means that $\Phi_{\sigma}$ is obtained from $\pi^{-1}(\wt\Phi_{\wt\sigma})$ by adding those roots
from $C_1\cup R_0$, which become negative under the action of involution $\sigma$.

In this case $\Phi_{\sigma}\cap C_1=C_1$ consists of $(2n-1)$ roots. Assume that $\alpha=\epsi_i\in R_0$, $2\leq
i\leq n$. Obviously, $r_{\beta_1}(\alpha)=\alpha$, so
$\sigma(\alpha)=\sigma_2(\alpha)$. But
$\sigma_2(\alpha)=\pm\epsi_l<0$ if and only if
$\beta=\epsi_i+\epsi_l\in\Supp{\sigma_2}$ for some $l>1$
(and each such $\beta$ adds to $\Phi_{\sigma}$ two roots, $\epsi_i$ and $\epsi_l$). Thus,
$$
\begin{array}{ll}
l(\sigma)&=|\Phi_{\sigma}|=|\pi^{-1}(\wt\Phi_{\wt\sigma})|+
|C_1| + 2\cdot\#\{\beta\in\Supp{\sigma_2}\mid\row(\beta)<0\}=\\
&=l(\wt\sigma)+(2n-1)+2d(\sigma)=l(\wt\sigma)+2(n+2d(\sigma))-1,
\end{array}
$$
as required.}

Now we are able to obtain a formula describing the dimension of an orbit associated with an orthogonal rook placement. \theop{Let
$\Phi=B_n$ or $D_n$\textup, $W$ be
its Weyl group\textup, $D=\{\beta_1\succ\ldots\succ\beta_t\}$ be an orthogonal rook placement $\Phi^+$\textup, $\sigma\in W$ be the involution with the support $\Supp{\sigma}=D$\textup, and
$\Omega=\Omega_{D,\xi}$ be the associated orbit \textup(for some $\xi$\textup).
Then the dimension of this orbit equals
$$
\dim\Omega=l(\sigma)-s(\sigma)-2d(\sigma).
$$
In particular\textup, $\dim\Omega$ does not depend on $\xi$.
\label{theo_dim_orbit}}{Let $\chi$ be the irreducible character of the group $U$ corresponding to the orbit $\Omega$ (see
(\ref{formula_char_orbit})). We will denote by $\deg\chi$ the
dimension of the corresponding representation, then
$\deg\chi=q^{\frac{1}{2}\dim\Omega}$ (see Theorem
\ref{theo_char_orbit}).

We will proceed by the induction on $n$ (the inductive base can be checked directly). In view of Theorem
\ref{theo_char_decompose} (see (\ref{formula_char_decompose})),
$$
\deg\chi=\deg\wt\chi\cdot[U:(U_1\rtimes V')].
$$
Here $\wt\Omega$ is the orbit associated with the orthogonal rook placement
$\wt D$ (and the corresponding $\wt\xi$), and $\wt\chi$ is the corresponding irreducible character of the group $\wt U$. By the inductive assumption,
$$\deg\wt\chi=l(\wt\sigma)-s(\wt\sigma)-2d(\wt\sigma).$$ We will consider the different places of the root $\beta_1$.

1. Let $\col(\beta_1)>1$. This is a trivial case: here
$l(\sigma)=l(\wt\sigma)$, $s(\sigma)=s(\wt\sigma)$,
$d(\sigma)=d(\wt\sigma)$, and $\dim\Omega=\dim\wt\Omega$.
Indeed, $\Phi_2^+=\Phi_1^+=\Phi^+\setminus C_1$; in view of
Corollary \ref{coro_dim_pol}, the dimension of $\Omega$ (respectively, of
$\wt\Omega$) is equal to the number of symbols $\pm$ in the diagram
$(\Phi,\Sb_{\sigma})$ (respectively, $(\wt\Phi,\wt\Sb)$, where
$\wt\Sb=\Sb_{\wt\sigma}$). But $\Sb{\mid}_{\Phi_2^+}=\wt\Sb\circ\pi$
(see the proof of Proposition \ref{prop_sub}), and
$\Sb(\alpha)=\bullet$ for any $\alpha\in C_1$, so
$$
\dim\wt\Omega=\#\{\alpha\in\Phi_2^+\mid\Sb(\alpha)=\pm\}
=\#\{\alpha\in\Phi^+\mid\Sb(\alpha)=\pm\}=\dim\Omega.
$$

2. Let $\col(\beta_1)=1,\row(\beta_1)=i$, $1\leq i\leq n$ (i.e.,
$\beta_1=\epsi_1-\epsi_i$). Here $d(\sigma)=d(\wt\sigma)$. It is clear that
$$
\begin{array}{ll}
&[U:(U_1\rtimes
V')]=q^{|S^-(\beta_1)|}=q^{i-2},\\
&\deg\chi=q^{\frac{1}{2}\dim\Omega}=\deg\wt\chi\cdot[U:(U_1\rtimes
V')]=q^{\frac{1}{2}\dim\wt\Omega+i-2}.
\end{array}
$$
We see that $\dim\Omega=\dim\wt\Omega+2(i-2)$.
On the other hand, using Lemma \ref{lemm_sigma_i_m_j} and the inductive assumption, we obtain
$$
\begin{array}{ll}
l(\sigma)-s(\sigma)-2d(\sigma)&=(l(\wt\sigma)+2i-3)-(s(\wt\sigma)+1)-2d(\wt\sigma)=\\
&=\dim\wt\Omega+2(i-2)=\dim\Omega.
\end{array}
$$

3. Let $\col(\beta_1)=1, \row(\beta_1)=-i$, $1\leq i\leq n$ (i.e.,
$\beta_1=\epsi_1+\epsi_i$). The line of reasoning is similar to the previous case. Namely, $d(\sigma)=d(\wt\sigma)$ and, taking into account Lemma
\ref{lemm_sigma_i_p_j}, we conclude that
$$
\begin{array}{ll}
&[U:(U_1\rtimes V')]=q^{|S^-(\beta_1)|}=q^{m-i-2},\\
&\deg\chi=q^{\frac{1}{2}\dim\Omega}=\deg\wt\chi\cdot[U:(U_1\rtimes
V')]=q^{\frac{1}{2}\dim\wt\Omega+m-i-2},\\
&\dim\Omega=\dim\wt\Omega+2(m-i-2),\\
&l(\sigma)-s(\sigma)-2d(\sigma)=(l(\wt\sigma)+2(m-i)-3)-(s(\wt\sigma)+1)-2d(\wt\sigma)=\\
&=\dim\wt\Omega+2(m-i-2)=\dim\Omega.
\end{array}
$$

4. Let $\col(\beta_1)=1,\row(\beta_1)=0$ (i.e., $\Phi=B_n$ and
$\beta_1=\epsi_1$). Then $d(\wt\sigma)=0$, therefore, taking into account Lemma \ref{lemm_sigma_i}, we see that
$$
\begin{array}{ll}
&[U:(U_1\rtimes V')]=q^{|S^-(\beta_1)|}=q^{n-1},\\
&\deg\chi=q^{\frac{1}{2}\dim\Omega}=\deg\wt\chi\cdot[U:(U_1\rtimes
V')]=q^{\frac{1}{2}\dim\wt\Omega+n-1},\\
&\dim\Omega=\dim\wt\Omega+2(n-1),\\
&l(\sigma)-s(\sigma)-2d(\sigma)=(l(\wt\sigma)+2(n+d(\sigma))-1)-(s(\wt\sigma)+1)-2d(\sigma)=\\
&=\dim\wt\Omega+2(n-1)=\dim\Omega.
\end{array}
$$
The proof is complete.}\newpage

\exam{Let $\Phi=B_6$, $\sigma$ be the involution from Example
\ref{exam_inv_B_6}:
$$\Supp{\sigma}=D=\{\epsi_1,~\epsi_2+\epsi_6,~\epsi_3+\epsi_5\}.$$
According to Corollary \ref{coro_dim_pol}, the dimension of any orbit 
associated with $D$ is equal to the number of symbols $\pm$ in the
diagram $(\Phi,\Sb_D)$, so $\dim\Omega=18$. At the same time,
$$
\begin{array}{ll}
\Phi_{\sigma}=C_1&\cup\{\epsi_2-\epsi_4,~\epsi_2,~\epsi_2+\epsi_3,~
\epsi_2+\epsi_5,~\epsi_2+\epsi_6\}\\&\cup\{\epsi_3-\epsi_4,~\epsi_3,~
\epsi_3+\epsi_5,~\epsi_3+\epsi_6\}\cup\{\epsi_4+\epsi_5,~\epsi_4+\epsi_6\}\cup
\{\epsi_5,~\epsi_5+\epsi_6\}\cup\{\epsi_6\}.
\end{array}
$$
Hence, $l(\sigma)=|\Phi_{\sigma}|=25$,
$s(\sigma)=|\Supp{\sigma}|=3$, therefore, by Theorem
\ref{theo_dim_orbit}, we get that
$$
\dim\Omega=l(\sigma)-s(\sigma)-2d(\sigma)=25-3-2\cdot2=18.
$$}

%The author expresses his sincere gratitude to his supervisor
%A.N. Panov for setting the task and attention to the work.

\begin{comment}

\end{comment}

%\vspace{5mm}
%\noindent Поступила в редакцию 24/{\it V}/2008;\\
%\noindent в окончательном варианте --- 16/{\it VI}/2008.
%\vspace{5mm}

%\begin{center}
%\maintitle[Communicated by Dr. Sci. (Phys. \& Math.) Prof.
%V.E.\,Voskresenskii] {BASIC SUBSYSTEMS OF ROOT SYSTEMS OF TYPES
%$B_n$ AND $D_n$\\AND ASSOCIATED COADJOINT ORBITS}
%\authorright[2008]
%\anauthor[Ignatev Mikhail Viktorovich
%\email{mihail\_ignatev@mail.ru}, Dept. of Algebra and Geometry,
%Samara State University, Samara, 443011, Russia.]{M.V.\,Ignatev}{}
%\end{center}

\bigskip\textsc{Mikhail V. Ignatev: National Research University Higher School of Economics, Pokrovsky Boulevard 11, 109028, Moscow, Russia}

\emph{E-mail address}: \texttt{mihail.ignatev@gmail.com}

\bigskip\textsc{Mikhail S. Venchakov: National Research University Higher School of Economics, Pokrovsky Boulevard 11, 109028, Moscow, Russia}

\emph{E-mail address}: \texttt{mihail.venchakov@gmail.com}

\end{document}